\documentclass[11pt]{amsart}
\usepackage[T1]{fontenc}
\usepackage{lmodern}
\usepackage[margin=1in]{geometry}
\usepackage{amsmath,amssymb,mathtools}
\usepackage{tikz-cd}
\usepackage{enumitem}
\usepackage{microtype}
\usepackage[hidelinks]{hyperref}
\hypersetup{
  pdftitle={Admissible balanced pairs over cleft extensions},
  pdfauthor={Junpeng Ren, Zheng-Gang Wu, Mengzhen Zhu},
  pdfsubject={Balanced pairs and cotorsion pairs in cleft extensions}
}

\newtheorem{theorem}{Theorem}[section]
\newtheorem{lemma}[theorem]{Lemma}
\newtheorem{proposition}[theorem]{Proposition}
\newtheorem{corollary}[theorem]{Corollary}
\theoremstyle{definition}
\newtheorem{definition}[theorem]{Definition}
\newtheorem{example}[theorem]{Example}
\theoremstyle{remark}
\newtheorem{remark}[theorem]{Remark}
\theoremstyle{plain}
\newtheorem*{theoremA}{Theorem A}
\newtheorem*{theoremB}{Theorem B}
\newtheorem*{theoremC}{Theorem C}
\newtheorem*{theoremD}{Theorem D}

\numberwithin{equation}{section}
\setlist[enumerate,1]{label=\textup{(\arabic*)},leftmargin=2.2em}
\title[Balanced pairs over cleft extensions]
{Admissible balanced pairs over cleft extensions}
\author[J. Ren]{Junpeng Ren}
\author[Z.-G. Wu]{Zheng-Gang Wu}
\author[M. Zhu]{Mengzhen Zhu$^{*}$}
\thanks{$^{*}$Corresponding author.}

\begin{document}

\begin{abstract}
	Let $(\mathcal{B}, \mathcal{A}, i, e, l)$ be a cleft extension of abelian categories such that
	$e$ has a right adjoint $r$. It is proved that $(\mathcal{C}, \mathcal{D})$ is an admissible
	balanced pair in $\mathcal{B}$ if and only if
	$(\operatorname{add}(l\mathcal{C}), \operatorname{add}(r\mathcal{D}))$ is an admissible
	balanced pair in $\mathcal{A}$. When the extension is directed with respect to a finite
	product decomposition of $\mathcal{B}$, the lifting of complete hereditary cotorsion pairs
	$(\mathcal{C}, \mathcal{Y})$ and $(\mathcal{X}, \mathcal{D})$ to $\mathcal{A}$ is
	characterized by the vanishing of first derived functors. Under a compatibility condition,
	the pair of special classes $(\Delta(\mathcal{C}), \nabla(\mathcal{D}))$ is then shown to be
	admissible balanced if and only if $(\mathcal{C}, \mathcal{D})$ is.
	As an application, for rings $A, B$ and a $(B, A)$-bimodule $U$ of finite flat dimension
	over $A$ and finite projective dimension over $B$, the triangular matrix ring
	$\left(\begin{smallmatrix}A & 0\\ U & B\end{smallmatrix}\right)$ is left virtually
	Gorenstein if and only if $A$ and $B$ are.
\end{abstract}

\maketitle

\section{Introduction}\label{sec:intro}

Balance in relative homological algebra was studied by Enochs and Jenda
\cite{EnochsJenda}. Let $\mathcal{C}$ and $\mathcal{D}$ be subcategories of an abelian
category $\mathcal{B}$ such that $\mathcal{C}$ is precovering and $\mathcal{D}$ is preenveloping.
Following Chen \cite{Chen}, the pair $(\mathcal{C}, \mathcal{D})$ is called \emph{balanced}
if every object of $\mathcal{B}$ admits a $\mathcal{C}$-resolution that becomes acyclic
under $\operatorname{Hom}_{\mathcal{B}}(-, D)$ for every $D \in \mathcal{D}$, and a $\mathcal{D}$-coresolution that
becomes acyclic under $\operatorname{Hom}_{\mathcal{B}}(C, -)$ for every $C \in \mathcal{C}$. The relative
derived functors of $\operatorname{Hom}$ can then be computed using either the
resolutions or the coresolutions \cite[Lemma~2.1]{Chen}. A balanced pair $(\mathcal{C}, \mathcal{D})$ is called \emph{admissible}
if every $\mathcal{C}$-precover is an epimorphism and every $\mathcal{D}$-preenvelope is a
monomorphism. When $\mathcal{B}$ has small Ext groups, the admissible balanced pairs in $\mathcal{B}$
correspond bijectively to the exact structures on $\mathcal{B}$ with enough
projective and injective objects \cite{WangLiHuang}.

Balanced pairs also arise from cotorsion theory: if $(\mathcal{C}, \mathcal{W}, \mathcal{D})$ is a
complete hereditary cotorsion triplet, then $(\mathcal{C}, \mathcal{D})$ is an admissible
balanced pair \cite[Theorem~4.1]{EJTX}, \cite[Proposition~2.6]{Chen}, and
the converse is studied in \cite{EPZ}. An example comes from Gorenstein
homological algebra. Let $\operatorname{GProj} R$ and $\operatorname{GInj} R$ denote the classes of
Gorenstein projective and Gorenstein injective left modules over a ring
$R$. The ring $R$ is called \emph{left virtually Gorenstein}
\cite{BeligiannisReiten,EPZ} if
$(\operatorname{GProj} R)^\perp = {}^\perp(\operatorname{GInj} R)$, the orthogonal classes being taken
with respect to $\operatorname{Ext}^1_R$. Estrada, Hu and
Trlifaj \cite{EHT} proved that $R$ is left virtually Gorenstein if and
only if $(\operatorname{GProj} R, \operatorname{GInj} R)$ is a balanced pair. Balanced pairs have
also been defined in triangulated categories \cite{FuHuZhangZhu} and in
extriangulated categories \cite{HeZhou,XuFu}.

The present article is concerned with the passage of admissible balanced
pairs between an abelian category and an extension of it. For a formal
triangular matrix ring
\[
 T = \begin{pmatrix}A & 0\\ U & B\end{pmatrix},
\]
and for classes $\mathcal{C}_1, \mathcal{D}_1$ of $A$-modules and
$\mathcal{C}_2, \mathcal{D}_2$ of $B$-modules, Mao proved that the
$T$-modules induced from $\mathcal{C}_1$ and $\mathcal{C}_2$ and those
coinduced from $\mathcal{D}_1$ and $\mathcal{D}_2$ form an admissible balanced
pair over $T$ if and only if $(\mathcal{C}_1, \mathcal{D}_1)$ and
$(\mathcal{C}_2, \mathcal{D}_2)$ are admissible balanced pairs over $A$ and
$B$ \cite[Theorem~2.2]{Mao}. A left $T$-module is a triple $(X, Y, \varphi)$
with $\varphi: U \otimes_A X \to Y$; let $\widetilde\varphi$ denote the
adjoint of $\varphi$. Mao also considered the $T$-modules for which $\varphi$
is a monomorphism with $X \in \mathcal{C}_1$ and
$\operatorname{Coker}\varphi \in \mathcal{C}_2$, and those for which
$\widetilde\varphi$ is an epimorphism with
$\operatorname{Ker}\widetilde\varphi \in \mathcal{D}_1$ and
$Y \in \mathcal{D}_2$. If these two classes form an admissible balanced pair
over $T$, then so do $(\mathcal{C}_1, \mathcal{D}_1)$ and
$(\mathcal{C}_2, \mathcal{D}_2)$ \cite[Proposition~2.4]{Mao}, and the
converse holds under cotorsion-pair hypotheses and $\operatorname{Tor}$ and
$\operatorname{Ext}$ vanishing conditions \cite[Theorem~2.5]{Mao}. For
Gorenstein projective and Gorenstein injective modules, this gives
conditions under which $T$ is left virtually Gorenstein if and only if $A$
and $B$ are \cite[Theorem~2.8]{Mao}.

The module category of $T$ is a cleft extension of
$A\text{-}\operatorname{Mod} \times B\text{-}\operatorname{Mod}$. Cleft extensions of abelian categories
were introduced by Beligiannis \cite{BelCleft} as a common generalization
of cleft extensions of rings and of the trivial extensions of abelian
categories studied by Fossum, Griffith and Reiten \cite{FGR}. A cleft
extension $(\mathcal{B}, \mathcal{A}, i, e, l)$ consists of abelian categories $\mathcal{A}, \mathcal{B}$ and
additive functors
\[
\begin{tikzcd}[column sep=large]
 \mathcal{B}\arrow[r,"i"]&
 \mathcal{A}\arrow[r,"e"']&
 \mathcal{B}\arrow[l,bend right=35,"l"']
\end{tikzcd}
\qquad
\begin{gathered}
 l \dashv e, \qquad ei \cong \mathrm{id}_{\mathcal{B}}, \\
 e\text{ is faithful and exact}.
\end{gathered}
\]
For rings, a pair of homomorphisms $R\xrightarrow{\alpha}S\xrightarrow{\beta}R$
with $\beta\alpha = \mathrm{id}_R$ gives a cleft extension of module categories;
trivial extensions, triangular matrix rings, tensor rings and Morita
rings with zero pairings arise in this way. By
\cite{BelCleft,Marmaridis}, every cleft extension is equivalent to an
\emph{$\eta$-extension} $\mathcal{B} \ltimes_\eta F$, where $F: \mathcal{B} \to \mathcal{B}$ is a right
exact functor and $\eta: F^2 \to F$ is an associative natural
transformation, which will be called the multiplication. The objects of
$\mathcal{B} \ltimes_\eta F$ are the pairs $(X, f)$ with $f: FX \to X$ and
$f \circ Ff = f \circ \eta_X$, and $e(X, f) = X$.

The relative homology of cleft extensions was developed in
\cite{BelRelative}. Since then, finitistic and other homological
dimensions \cite{GreenPsaroudakisSolberg,MaZhengLiu}, Gorenstein objects
and singularity categories \cite{Kostas,LiangMaYang,Qin}, cotorsion pairs
\cite{ChenQinYang,Hu,LiuGao}, and silting and tilting objects
\cite{ZhaoSunSilting} have been compared between a cleft extension and
its base category. The purpose of the present article is to carry out
this comparison for admissible balanced pairs.

Suppose that $e$ has a right adjoint $r$. Then $F$ has a right adjoint
$G$, and $i$ has a left adjoint $q$ and a right adjoint $s$:
\[
 l \dashv e \dashv r, \qquad q \dashv i \dashv s.
\]
For the ring $T$, the functors $l$ and $r$ give the induced and coinduced
modules, and $q$ and $s$ will define the special classes.

Throughout the paper, subcategories are full, additive and closed under
isomorphisms and direct summands. For a class $\mathcal{U}$ of objects,
$\operatorname{add}\mathcal{U}$ denotes its closure under finite direct sums and direct
summands, and for a functor $\Gamma$, $\Gamma\mathcal{U}$ denotes the class of objects
isomorphic to $\Gamma U$ with $U \in \mathcal{U}$. The first main result of the paper may then be stated as follows.

\begin{theoremA}
Let $(\mathcal{B}, \mathcal{A}, i, e, l)$ be a cleft extension, and suppose that $e$ has a right
adjoint $r$. Then $(\mathcal{C}, \mathcal{D})$ is an admissible balanced pair in $\mathcal{B}$ if and
only if $(\operatorname{add}(l\mathcal{C}), \operatorname{add}(r\mathcal{D}))$ is an admissible balanced pair in $\mathcal{A}$.
\end{theoremA}

Theorem A requires no homological condition. It applies to every cleft
extension of rings, in particular to trivial extensions and to Morita
rings with zero pairings (Corollary~\ref{cor:split-rings}), and for comma
categories it recovers \cite[Theorem~2.2]{Mao}
(Corollary~\ref{cor:comma-induced}). The closure $\operatorname{add}$ is taken because
$l\mathcal{C}$ and $r\mathcal{D}$ need not be closed under direct summands
(Example~\ref{ex:summands}). On the coinduced side, $r$ cannot be replaced by
$i$: the pair $(\operatorname{add}(l\mathcal{C}), i\mathcal{D})$ may fail to be admissible balanced
even when $(\mathcal{C}, \mathcal{D})$ is (Example~\ref{ex:section}).

Suppose now that $\mathcal{B}$ has enough projective and
injective objects, and put
\[
 \Delta(\mathcal{C}) = \{M \in \mathcal{A} \mid L_1q(M) = 0, \ qM \in \mathcal{C}\}, \qquad
 \nabla(\mathcal{D}) = \{N \in \mathcal{A} \mid R^1s(N) = 0, \ sN \in \mathcal{D}\}.
\]
Following Hu
\cite[Subsections~3.1 and~3.6]{Hu}, they are called \emph{special classes}. In the model $\mathcal{B} \ltimes_\eta F$, the
condition $L_1q(X, f) = 0$ means that
\[
 F^2X\xrightarrow{Ff-\eta_X}FX\xrightarrow{f}X
\]
is exact at $FX$ (Proposition~\ref{prop:low-degree}); for the ring $T$,
the structure map is then a monomorphism. One has
$l\mathcal{C} \subseteq \Delta(\mathcal{C})$ and $r\mathcal{D} \subseteq \nabla(\mathcal{D})$, and both
inclusions can be strict (Example~\ref{ex:path}).

If $(\Delta(\mathcal{C}), \nabla(\mathcal{D}))$ is an admissible balanced pair, then so is
$(\mathcal{C}, \mathcal{D})$ (Lemma~\ref{lem:reflection}). For the converse, let
$(\mathcal{C}, \mathcal{Y})$ and $(\mathcal{X}, \mathcal{D})$ be complete hereditary cotorsion pairs in
$\mathcal{B}$, and write $e^{-1}(\mathcal{U}) = \{M \in \mathcal{A} \mid eM \in \mathcal{U}\}$ for a class
$\mathcal{U} \subseteq \mathcal{B}$. The problem is to lift them to complete hereditary
cotorsion pairs
\[
 (\Delta(\mathcal{C}), e^{-1}(\mathcal{Y})), \qquad (e^{-1}(\mathcal{X}), \nabla(\mathcal{D}))
\]
in $\mathcal{A}$. For $\eta$-extensions, the pair
$({}^\perp e^{-1}(\mathcal{Y}), e^{-1}(\mathcal{Y}))$ is known to be complete when
$L_1F(\mathcal{C}) = 0$ and $F(\mathcal{C}) \subseteq \mathcal{Y}$
\cite[Theorem~3.23]{ChenQinYang}, and earlier for right trivial
extensions \cite[Theorem~3.13]{Hu} and image-nilpotent extensions
\cite[Theorem~4.13]{LiuGao}. For the arrow category of $\mathcal{B}$ with the
same cotorsion pair at both vertices, the second condition reads
$\mathcal{C} \subseteq \mathcal{Y}$. For comma categories the lifting is known without it
\cite{HuZhuComma,YuanHeWu}. Some condition on $F$ is nevertheless needed:
by Example~\ref{ex:nilpotent-obstruction}, the identification
${}^\perp e^{-1}(\mathcal{Y}) = \Delta(\mathcal{C})$, on which the lifting rests, can fail
even when $F$ is exact and $F^2 = 0$.

Suppose that $\mathcal{B} = \prod_{j = 1}^n\mathcal{B}_j$, and let $F_{ij}: \mathcal{B}_j \to \mathcal{B}_i$ be
the matrix components of $F$. The cleft extension $\mathcal{A} = \mathcal{B} \ltimes_\eta F$
is called \emph{finite directed} if $F_{ij} = 0$ for all $i \leq j$; no
condition is imposed on the multiplication $\eta$. For $n = 2$ these are
the comma categories, and in particular the module categories of formal
triangular matrix rings; related constructions for quiver
representations and linear and generalized linear Reedy categories have been studied in
\cite{DaleziosStovicek,DiLiLiang,HolmJorgensen}.

For such extensions, every object
of $e^{-1}(\mathcal{Y})$ has a finite filtration with factors in $i\mathcal{Y}$, so that
$\Delta(\mathcal{C}) = {}^\perp e^{-1}(\mathcal{Y})$ without any vanishing condition
(Proposition~\ref{prop:directed-pair}). The lifting problem then has the
following answer.

\begin{theoremB}
Let $\mathcal{A} = \mathcal{B} \ltimes_\eta F$ be a finite directed cleft extension such that
$e$ has a right adjoint, with $\mathcal{B} = \prod_{j = 1}^n\mathcal{B}_j$, and let $G$ be the
right adjoint of $F$. Suppose that each $\mathcal{B}_j$ has
enough projective and injective objects. For each $j$, let
$(\mathcal{C}_j, \mathcal{Y}_j)$ and $(\mathcal{X}_j, \mathcal{D}_j)$ be complete hereditary cotorsion pairs
in $\mathcal{B}_j$, and let $\mathcal{C}, \mathcal{Y}, \mathcal{X}, \mathcal{D}$ be the corresponding product classes.
\begin{enumerate}
 \item $(\Delta(\mathcal{C}), e^{-1}(\mathcal{Y}))$ is a complete hereditary cotorsion pair
 in $\mathcal{A}$ if and only if $L_1F(\mathcal{C}) = 0$.
 \item $(e^{-1}(\mathcal{X}), \nabla(\mathcal{D}))$ is a complete hereditary cotorsion pair
 in $\mathcal{A}$ if and only if $R^1G(\mathcal{D}) = 0$.
\end{enumerate}
\end{theoremB}

The main result on the special classes may then be stated as follows.

\begin{theoremC}
In the setting of Theorem B, suppose that
$\mathcal{C}_j \cap \mathcal{Y}_j \subseteq \mathcal{X}_j$ and $\mathcal{X}_j \cap \mathcal{D}_j \subseteq \mathcal{Y}_j$
for each $j$. If
\[
 L_1F(\mathcal{C}) = 0, \qquad R^1G(\mathcal{D}) = 0,
\]
then the following are equivalent:
\begin{enumerate}
 \item $(\mathcal{C}, \mathcal{D})$ is admissible balanced in $\mathcal{B}$;
 \item $(\operatorname{add}(l\mathcal{C}), \operatorname{add}(r\mathcal{D}))$ is admissible balanced in $\mathcal{A}$;
 \item $(\Delta(\mathcal{C}), \nabla(\mathcal{D}))$ is admissible balanced in $\mathcal{A}$;
 \item $\mathcal{Y} = \mathcal{X}$;
 \item $\Delta(\mathcal{C})^\perp = {}^\perp\nabla(\mathcal{D})$.
\end{enumerate}
\end{theoremC}

The equivalence of (1) and (2) is Theorem A, and the conditions
$\mathcal{C}_j \cap \mathcal{Y}_j \subseteq \mathcal{X}_j$ and $\mathcal{X}_j \cap \mathcal{D}_j \subseteq \mathcal{Y}_j$ are
only used in the proof of $(1) \Rightarrow (4)$. Condition (5) restates (4) in
$\mathcal{A}$, since $\Delta(\mathcal{C})^\perp = e^{-1}(\mathcal{Y})$ and
${}^\perp\nabla(\mathcal{D}) = e^{-1}(\mathcal{X})$ by Theorem B. When $\mathcal{Y} = \mathcal{X}$, the
lifted pairs form a complete hereditary cotorsion triplet, and a
criterion of \cite{EPZ} connects this equality with balance at the three
levels. For triangular matrix rings Theorem C recovers \cite[Theorem~2.5]{Mao}.

Suppose now that $U_A$ has finite flat dimension and ${}_BU$ has finite
projective dimension. For Artin algebras, Eshraghi, Hafezi, Salarian and Li
\cite[Proposition~3.3]{EHSL} proved that $T$ is virtually Gorenstein if and
only if $A$ and $B$ are; their result is stated under four compatibility
conditions on $U$, which hold in this situation.
Mao \cite[Theorem~2.8]{Mao} characterized virtually Gorenstein triangular
matrix rings when $U_A$ is flat, ${}_BU$ is projective, and $A$ and $B$
are left noetherian rings of finite left self-injective dimension.

For arbitrary rings $A$ and $B$, the Gorenstein projective and Gorenstein injective
$T$-modules form the special classes
$\Delta(\operatorname{GProj} A \times \operatorname{GProj} B)$ and
$\nabla(\operatorname{GInj} A \times \operatorname{GInj} B)$ \cite{LiZhengHuZhu,YuanHeWu}, and
the vanishing conditions hold by dimension shifting. Together with the
characterization of virtually Gorenstein rings in \cite{EHT}, the lifting
of cotorsion triplets then yields the following result.

\begin{theoremD}
Let $A$ and $B$ be rings and $U$ a $(B, A)$-bimodule such that $U_A$ has
finite flat dimension and ${}_BU$ has finite projective dimension. Put
$T = \left(\begin{smallmatrix}A & 0\\ U & B\end{smallmatrix}\right)$. Then $T$
is left virtually Gorenstein if and only if $A$ and $B$ are left
virtually Gorenstein. Equivalently, $(\operatorname{GProj} T, \operatorname{GInj} T)$ is admissible
balanced if and only if $(\operatorname{GProj} A, \operatorname{GInj} A)$ and $(\operatorname{GProj} B, \operatorname{GInj} B)$ are
admissible balanced.
\end{theoremD}

In particular, the ring of lower triangular $n \times n$ matrices over a
ring $R$ is left virtually Gorenstein if and only if $R$ is
(Corollary~\ref{cor:lower-triangular}).

The paper is organized as follows. Section~\ref{sec:preliminaries}
recalls cleft extensions and balanced pairs, and contains the low-degree
exact sequence on which the special classes are based.
Section~\ref{sec:induced} proves Theorem A. Section~\ref{sec:structural}
treats the special classes and proves Theorems B and C. Section~\ref{sec:applications} contains the applications to
cleft extensions of rings, comma categories and triangular matrix rings,
and the proof of Theorem D.

\section{Preliminaries}\label{sec:preliminaries}

Throughout the paper we keep the conventions on subcategories and the
notation $\operatorname{add}\mathcal{U}$ and $\Gamma\mathcal{U}$ introduced in
Section~\ref{sec:intro}.

\subsection{Cleft extensions}\label{subsec:cleft}

Let $\mathcal{A}$ and $\mathcal{B}$ be abelian categories. Following
\cite{BelCleft,BelRelative,Marmaridis}, a \emph{cleft extension}
$(\mathcal{B}, \mathcal{A}, i, e, l)$ consists of additive functors
\[
 \mathcal{B}\xrightarrow{i}\mathcal{A}\xrightarrow{e}\mathcal{B}, \qquad l: \mathcal{B} \to \mathcal{A},
 \qquad
 l \dashv e, \quad ei \cong \mathrm{id}_{\mathcal{B}},
\]
where $e$ is faithful and exact. The functor $i$ is then fully faithful and exact, and has a left adjoint $q$. 
Throughout, we assume that $e$ admits a right adjoint $r$. Then $(\mathcal{B},\mathcal{A},i,e,r)$ is a cleft coextension, and hence $i$ admits a right adjoint $s$; see \cite[Section~2.4]{BelRelative}.
These functors fit into the diagram
\begin{equation}\label{eq:adjunctions}
\begin{tikzcd}[column sep=large]
 \mathcal{B}\arrow[r,"i" description]&
 \mathcal{A}\arrow[l,bend right=40,"q"']
    \arrow[l,bend left=40,"s"]
    \arrow[r,"e" description]&
 \mathcal{B}\arrow[l,bend right=40,"l"']
    \arrow[l,bend left=40,"r"]
\end{tikzcd}
\qquad
\begin{gathered}
 q \dashv i \dashv s, \\
 l \dashv e \dashv r.
\end{gathered}
\end{equation}
Naturally isomorphic functors are identified when this causes no ambiguity.

Let $\nu: \mathrm{id}_{\mathcal{B}} \to el$ and $\mu: le \to \mathrm{id}_{\mathcal{A}}$ be the unit and counit
of $l \dashv e$. Since $e(\mu_M)\nu_{eM} = \mathrm{id}_{eM}$, each $\mu_M$ is an
epimorphism. Put $\Omega(M) = \operatorname{Ker}\mu_M$ and $F = e\Omega i$. Applying $e$
to the exact sequence $0 \to \Omega \to le\xrightarrow{\mu}\mathrm{id}_{\mathcal{A}} \to 0$
at objects $iX$ yields a split exact sequence
\[
 0 \longrightarrow F \longrightarrow el
   \xrightarrow{e\mu i}\mathrm{id}_{\mathcal{B}} \longrightarrow 0,
 \qquad el \cong \mathrm{id}_{\mathcal{B}} \oplus F,
\]
so that $F$ is additive and right exact. Since $F$ is a direct summand of $el$ and $el$ has the
right adjoint $er$, the functor $F$ has a right adjoint $G$. Since $el \dashv er$, uniqueness of adjoints
yields $er \cong \mathrm{id}_{\mathcal{B}} \oplus G$.

\begin{lemma}\label{lem:adjoints}
The functors in \eqref{eq:adjunctions} satisfy
\[
 ql \cong \mathrm{id}_{\mathcal{B}}, \quad sr \cong \mathrm{id}_{\mathcal{B}}, \quad
 qi \cong \mathrm{id}_{\mathcal{B}}, \quad si \cong \mathrm{id}_{\mathcal{B}}.
\]
The functors $l$ and $q$ preserve projective objects, and $r$ and $s$
preserve injective objects. If $\mathcal{B}$ has enough projective objects, then
so does $\mathcal{A}$. Dually, if $\mathcal{B}$ has enough injective objects, then so does
$\mathcal{A}$.
\end{lemma}

These are standard properties of cleft extensions and coextensions;
see \cite[Propositions~2.3 and~5.4(1)]{BelCleft},
\cite[Section~2]{Kostas}, and their duals.

By monadicity of $e$, every cleft extension is equivalent to an
$\eta$-extension $\mathcal{B} \ltimes_\eta F$ \cite{BelCleft,Marmaridis}. Here
$\eta: F^2 \to F$ is an associative natural transformation, that is,
$\eta \circ F\eta = \eta \circ \eta F$. An object of $\mathcal{B} \ltimes_\eta F$ is a
pair $(X, f)$ with $f: FX \to X$, and a morphism
$a: (X, f) \to (Y, h)$ is a morphism $a: X \to Y$ in $\mathcal{B}$. These satisfy,
respectively,
\[
 f \circ Ff = f \circ \eta_X, \qquad af = h \circ Fa.
\]
In this model, four of the functors in \eqref{eq:adjunctions} are given by
\begin{equation}\label{eq:model}
 \begin{alignedat}{2}
 e(X, f) & = X, & \qquad q(X, f) & = \operatorname{Coker} f, \\[2pt]
 iX & = (X, 0), & lX & = \left(X \oplus FX,
       \begin{pmatrix}0 & 0\\ \mathrm{id}_{FX} & \eta_X\end{pmatrix}\right).
 \end{alignedat}
\end{equation}
Dually, let $\zeta: G \to G^2$ be the mate of $\eta$. Under the
correspondence between morphisms $f: FX \to X$ and $g: X \to GX$, the category
$\mathcal{B} \ltimes_\eta F$ is equivalent to the coextension $G \rtimes_\zeta\mathcal{B}$,
whose objects are pairs $[X, g]$ with $Gg \circ g = \zeta_X \circ g$. In this
model,
\[
 s[X, g] = \operatorname{Ker} g, \qquad
 rX = \left[X \oplus GX,
       \begin{pmatrix}0 & \mathrm{id}_{GX}\\0 & \zeta_X\end{pmatrix}\right];
\]
see \cite[Sections~2.2 and~2.4]{BelRelative}.

\subsection{Low-degree comparison}

The special classes in Section~\ref{sec:structural} involve the
vanishing of $L_1q$ and $R^1s$. To relate these conditions to
$\operatorname{Ext}^1$-orthogonality we use a Grothendieck spectral sequence, which is
the absolute case of the relative spectral sequence of Beligiannis
\cite[Proposition~7.1]{BelRelative}.

\begin{proposition}\label{prop:low-degree}
Let $(\mathcal{B},\mathcal{A},i,e,l)$ be a cleft extension of abelian categories.
Suppose that $\mathcal{B}$ has enough projective and enough injective objects. For every $M\in\mathcal{A}$ and $Y\in\mathcal{B}$, there is an exact sequence, natural in both variables,
\begin{align}
 0 & \longrightarrow \operatorname{Ext}^1_{\mathcal{B}}(qM, Y)
\longrightarrow \operatorname{Ext}^1_{\mathcal{A}}(M, iY)
\longrightarrow \operatorname{Hom}_{\mathcal{B}}(L_1q(M), Y) \notag\\
 & \longrightarrow \operatorname{Ext}^2_{\mathcal{B}}(qM, Y)
\longrightarrow \operatorname{Ext}^2_{\mathcal{A}}(M, iY).
 \label{eq:five-term}
\end{align}
In particular, for every injective $I \in \mathcal{B}$, there is a natural isomorphism
\begin{equation}\label{eq:injective-test}
 \operatorname{Ext}^1_{\mathcal{A}}(M, iI) \cong \operatorname{Hom}_{\mathcal{B}}(L_1q(M), I).
\end{equation}
In the model \eqref{eq:model},
\begin{equation}\label{eq:L1q-model}
 L_1q(X, f) \cong
 \operatorname{Ker} f/\operatorname{Im}(Ff-\eta_X).
\end{equation}

\end{proposition}

\begin{proof}
Let $h_Y = \operatorname{Hom}_{\mathcal{B}}(-, Y)$. The adjunction $q \dashv i$ gives a natural
isomorphism $h_Yq \cong \operatorname{Hom}_{\mathcal{A}}(-, iY)$. By Lemma~\ref{lem:adjoints}, the
category $\mathcal{A}$ has enough projective objects and $q$ preserves projective
objects, so that $\operatorname{Ext}^n_{\mathcal{B}}(qP, Y) = 0$ for every projective object $P$
of $\mathcal{A}$ and every $n \geq 1$; that is, $qP$ is $h_Y$-acyclic. We therefore
have a Grothendieck spectral sequence
\[
 E_2^{p, t} = \operatorname{Ext}^p_{\mathcal{B}}(L_tq(M), Y)
 \quad\Longrightarrow\quad
 \operatorname{Ext}^{p+t}_{\mathcal{A}}(M, iY).
\]
Its five-term exact sequence is
\[
 0 \longrightarrow E_2^{1, 0}
\longrightarrow \operatorname{Ext}^1_{\mathcal{A}}(M, iY)
\longrightarrow E_2^{0, 1}
\longrightarrow E_2^{2, 0}
\longrightarrow \operatorname{Ext}^2_{\mathcal{A}}(M, iY),
\]
and substituting the description of the $E_2$-terms yields
\eqref{eq:five-term}. Now let $I$ be injective in $\mathcal{B}$. Then
$\operatorname{Ext}^n_{\mathcal{B}}(-, I) = 0$ for all $n \geq 1$, so that \eqref{eq:five-term}
with $Y = I$ gives the natural isomorphism \eqref{eq:injective-test}.
Finally, \eqref{eq:L1q-model} follows from
\cite[Corollary~5.10(i)]{BelCleft}.

\end{proof}

Dually, for cleft coextensions of abelian categories, we obtain the following result.
\begin{proposition}\label{prop:low-degree-dual}
	Let $(\mathcal{B},\mathcal{A},i,e,r)$ be a cleft coextension of abelian categories.
	Suppose that $\mathcal{B}$ has enough projective and enough injective objects. For every $M\in\mathcal{A}$ and $Y\in\mathcal{B}$, there is an exact sequence, natural in both variables,
	\begin{align}
		0 & \longrightarrow \operatorname{Ext}^1_{\mathcal{B}}(Y, sM)
		\longrightarrow \operatorname{Ext}^1_{\mathcal{A}}(iY, M)
		\longrightarrow \operatorname{Hom}_{\mathcal{B}}(Y, R^1s(M)) \notag\\
		& \longrightarrow \operatorname{Ext}^2_{\mathcal{B}}(Y, sM)
		\longrightarrow \operatorname{Ext}^2_{\mathcal{A}}(iY, M).
		\label{eq:five-term-dual}
	\end{align}
	In particular, for every projective $P \in \mathcal{B}$, there is a natural isomorphism
	\begin{equation}\label{eq:projective-test}
		\operatorname{Ext}^1_{\mathcal{A}}(iP, M) \cong \operatorname{Hom}_{\mathcal{B}}(P, R^1s(M)).
	\end{equation}
	In the dual model of \eqref{eq:model},
	\begin{equation}\label{eq:R1s-model}
		R^1s[X, g] \cong
		\operatorname{Ker}(Gg-\zeta_X)/\operatorname{Im} g.
	\end{equation}
\end{proposition}

\subsection{Balanced pairs and cotorsion triplets}

Let $\mathcal{U}$ be a subcategory of an abelian category $\mathcal{B}$.
A morphism $u: U \to B$ with $U \in \mathcal{U}$ is a
\emph{$\mathcal{U}$-precover} if $\operatorname{Hom}_{\mathcal{B}}(U', u)$ is surjective for every
$U' \in \mathcal{U}$; it is \emph{epic} if $u$ is an epimorphism. Dually, a
morphism $v: B \to U$ with $U \in \mathcal{U}$ is a
\emph{$\mathcal{U}$-preenvelope} if $\operatorname{Hom}_{\mathcal{B}}(v, U')$ is surjective for
every $U' \in \mathcal{U}$; it is \emph{monic} if $v$ is a monomorphism.
The subcategory $\mathcal{U}$ is \emph{epic precovering} if every object
admits an epic $\mathcal{U}$-precover; \emph{monic preenveloping} is
defined dually.

Following Chen \cite{Chen}, a pair $(\mathcal{C}, \mathcal{D})$ of subcategories is an
\emph{admissible balanced pair} if $\mathcal{C}$ is epic precovering, $\mathcal{D}$ is monic
preenveloping, and every $B \in \mathcal{B}$ admits exact complexes
\[
 \cdots \longrightarrow C_1 \longrightarrow C_0 \longrightarrow B
\longrightarrow 0, \qquad
 0 \longrightarrow B \longrightarrow D^0 \longrightarrow D^1
\longrightarrow \cdots,
\]
with $C_j \in \mathcal{C}$ and $D^j \in \mathcal{D}$, both of which remain exact under
$\operatorname{Hom}_{\mathcal{B}}(\mathcal{C}, -)$ and $\operatorname{Hom}_{\mathcal{B}}(-, \mathcal{D})$. The following
characterization is used throughout in place of the definition.

\begin{lemma}[{\cite[Proposition~2.2]{Chen}, \cite[Lemma~3.1]{EPZ}}]\label{lem:balance-criterion}
Let $\mathcal{C}$ be epic precovering and $\mathcal{D}$ monic preenveloping. Then
$(\mathcal{C}, \mathcal{D})$ is admissible balanced if and only if, for every short exact
sequence $\sigma$,
\begin{equation}\label{eq:balance-criterion}
 \operatorname{Hom}_{\mathcal{B}}(\mathcal{C}, \sigma)\text{ is exact}
 \quad \Longleftrightarrow \quad
 \operatorname{Hom}_{\mathcal{B}}(\sigma, \mathcal{D})\text{ is exact}.
\end{equation}
\end{lemma}

For an admissible balanced pair $(\mathcal{C}, \mathcal{D})$ in $\mathcal{B}$, let $\mathcal{E}_{\mathcal{B}}$ denote
the common class of short exact sequences in \eqref{eq:balance-criterion}.
This is the exact structure on $\mathcal{B}$ projectively generated by $\mathcal{C}$.
When $\mathcal{B}$ has small Ext groups, Wang, Li and
Huang \cite[Theorem~2.2]{WangLiHuang} showed that the assignment
$(\mathcal{C}, \mathcal{D}) \mapsto \mathcal{E}_{\mathcal{B}}$ gives a bijection between admissible balanced
pairs and exact structures on $\mathcal{B}$ with enough projective and injective
objects.

For a class $\mathcal{U}$ of objects, write
\[
 {}^\perp\mathcal{U} = \{B \mid \operatorname{Ext}^1_{\mathcal{B}}(B, U) = 0\text{ for all }U \in \mathcal{U}\},
 \qquad
 \mathcal{U}^\perp = \{B \mid \operatorname{Ext}^1_{\mathcal{B}}(U, B) = 0\text{ for all }U \in \mathcal{U}\}.
\]
A \emph{cotorsion pair} $(\mathcal{C}, \mathcal{Y})$ \cite[Chapter~7]{EnochsJenda} satisfies $\mathcal{C} = {}^\perp\mathcal{Y}$ and
$\mathcal{Y} = \mathcal{C}^\perp$. It is \emph{complete} if every object $B$ occurs in short
exact sequences
\[
 0 \to Y \to C \to B \to 0, \qquad
 0 \to B \to Y' \to C' \to 0
\]
with $C, C' \in \mathcal{C}$ and $Y, Y' \in \mathcal{Y}$; these are called special precover and
special preenvelope sequences. It is \emph{hereditary} if
$\operatorname{Ext}^m_{\mathcal{B}}(\mathcal{C}, \mathcal{Y}) = 0$ for all $m \geq 1$. A triple $(\mathcal{C}, \mathcal{W}, \mathcal{D})$ is a
\emph{complete hereditary cotorsion triplet} if both $(\mathcal{C}, \mathcal{W})$ and
$(\mathcal{W}, \mathcal{D})$ are complete hereditary cotorsion pairs.

The following result relates complete hereditary cotorsion triplets to admissible balanced pairs; see \cite[Theorem~4.1]{EJTX} and \cite[Proposition~2.6]{Chen}.

\begin{lemma}\label{lem:triplet-balanced}
Suppose that $\mathcal{B}$ has enough projective and injective objects. If
$(\mathcal{C}, \mathcal{W}, \mathcal{D})$ is a complete hereditary cotorsion triplet in $\mathcal{B}$, then
$(\mathcal{C}, \mathcal{D})$ is an admissible balanced pair in $\mathcal{B}$.
\end{lemma}

\section{Admissible balanced pairs under cleft extensions}\label{sec:induced}
In this section, we study the transfer of admissible balanced pairs along a cleft extension of abelian categories.
Let $(\mathcal{B}, \mathcal{A}, i, e, l)$ be a cleft extension of abelian categories and suppose that $e$ admits a right adjoint $r$. Based on the discussion in Subsection~\ref{subsec:cleft}, there are two adjoint triples $q \dashv i \dashv s$ and $l \dashv e \dashv r$.
\begin{lemma}\label{lem:transfer}
	Let $\varepsilon:qi\to\mathrm{id}_{\mathcal{B}}$ be the counit
	of $q\dashv i$, and let
	$\vartheta:\mathrm{id}_{\mathcal{B}}\to si$ be the unit of
	$i\dashv s$.
	
	Let $\mathcal{F}\subseteq\mathcal{A}$ and
	$\mathcal{C}\subseteq\mathcal{B}$ be subcategories satisfying
	$l\mathcal{C}\subseteq\mathcal{F}$ and
	$q\mathcal{F}\subseteq\mathcal{C}$.
	If $p:M_0\to iB$ is an $\mathcal{F}$-precover, then
	\[
	\bar p=\varepsilon_B\circ q(p):qM_0\longrightarrow B
	\]
	is a $\mathcal{C}$-precover. Moreover, $\bar p$ is an epimorphism whenever $p$ is.
	
	Dually, let $\mathcal{L}\subseteq\mathcal{A}$ and $\mathcal{D}\subseteq\mathcal{B}$ be subcategories satisfying
	$r\mathcal{D}\subseteq\mathcal{L}$ and
	$s\mathcal{L}\subseteq\mathcal{D}$.
	If $v:iB\to N_0$ is an $\mathcal{L}$-preenvelope, then
	\[
	\bar v=s(v)\circ\vartheta_B:B\longrightarrow sN_0
	\]
	is a $\mathcal{D}$-preenvelope. Moreover, $\bar v$ is a monomorphism whenever $v$ is.
\end{lemma}

\begin{proof}
	By Lemma~\ref{lem:adjoints} and the hypotheses,
	\[
	\mathcal{C}=ql\mathcal{C}\subseteq q\mathcal{F}
	\subseteq\mathcal{C},
	\qquad
	\mathcal{D}=sr\mathcal{D}\subseteq s\mathcal{L}
	\subseteq\mathcal{D}.
	\]
	Thus $q\mathcal{F}=\mathcal{C}$ and
	$s\mathcal{L}=\mathcal{D}$.
	
	In particular, $qM_0\in\mathcal{C}$. Let $b:qF\to B$ be a morphism, where $F\in\mathcal{F}$. Under the adjunction $q\dashv i$, it corresponds to a morphism $\widetilde b:F\to iB$ satisfying
	\[
	b=\varepsilon_B\circ q(\widetilde b).
	\]
	Since $p$ is an $\mathcal{F}$-precover, we can write
	$\widetilde b=p\circ c$ for some $c:F\to M_0$. Hence
	\[
	b=\varepsilon_B\circ q(p)\circ q(c)
	=\bar p\circ q(c).
	\]
	Therefore $\bar p$ is a $\mathcal{C}$-precover.
	
	Dually, $sN_0\in\mathcal{D}$, and it suffices to consider a morphism $a:B\to sL$ with $L\in\mathcal{L}$. Under the adjunction $i\dashv s$, it corresponds to a morphism $\widetilde a:iB\to L$ satisfying
	\[
	a=s(\widetilde a)\circ\vartheta_B.
	\]
	Since $v$ is an $\mathcal{L}$-preenvelope, we can write $\widetilde a=c'\circ v$ for some $c':N_0\to L$. Thus
	\[
	a=s(c')\circ s(v)\circ\vartheta_B
	=s(c')\circ\bar v.
	\]
	Therefore $\bar v$ is a $\mathcal{D}$-preenvelope.
	
	Finally, since $i$ is fully faithful, $\varepsilon_B$ and $\vartheta_B$ are isomorphisms. The functor $q$, being a left adjoint, preserves epimorphisms, whereas
	$s$, being a right adjoint, preserves monomorphisms. This proves the remaining assertions.
\end{proof}

\begin{theorem}\label{thm:induced}
	Let $(\mathcal{B},\mathcal{A},i,e,l)$ be a cleft extension of abelian categories, and suppose that $e\dashv r$.
	For subcategories $\mathcal{C},\mathcal{D}\subseteq\mathcal{B}$,
	the following are equivalent:
	\begin{enumerate}
		\item $(\mathcal{C},\mathcal{D})$ is an admissible balanced
		pair in $\mathcal{B}$;
		\item $(\operatorname{add}(l\mathcal{C}),
		\operatorname{add}(r\mathcal{D}))$
		is an admissible balanced pair in $\mathcal{A}$.
	\end{enumerate}
\end{theorem}

\begin{proof}
	Put
	\[
	\mathcal{F}=\operatorname{add}(l\mathcal{C}),
	\qquad
	\mathcal{L}=\operatorname{add}(r\mathcal{D}).
	\]
	
	Suppose first that $\mathcal{C}$ is epic precovering. For $M\in\mathcal{A}$, choose an epic $\mathcal{C}$-precover $a:C_0\to eM$. Under the adjunction $l\dashv e$, its adjoint
	\[
	\widetilde a=\mu_M\circ l(a):lC_0\to M
	\]
	is an $l\mathcal{C}$-precover, and hence an
	$\mathcal{F}$-precover. Here $\mu$ denotes the counit of $l\dashv e$. If $\nu$ denotes the corresponding unit, then
	$a=e(\widetilde a)\circ\nu_{C_0}$.
	Since $a$ is an epimorphism, so is $e(\widetilde a)$. The faithfulness of $e$ implies that $\widetilde a$ is an epimorphism. Thus $\mathcal{F}$ is epic precovering.
	
	Conversely, suppose that $\mathcal{F}$ is epic precovering. Since $q$ is additive and $ql\cong\mathrm{id}_{\mathcal{B}}$, we have
	\[
	l\mathcal{C}\subseteq\mathcal{F},
	\qquad
	q\mathcal{F}\subseteq
	\operatorname{add}(ql\mathcal{C})=\mathcal{C}.
	\]
	For each $B\in\mathcal{B}$, take an epic
	$\mathcal{F}$-precover of $iB$. By
	Lemma~\ref{lem:transfer}, it induces an epic
	$\mathcal{C}$-precover of $B$. Therefore
	\begin{equation}\label{eq:precover-transfer}
		\mathcal{C}\text{ is epic precovering}
		\quad\Longleftrightarrow\quad
		\mathcal{F}\text{ is epic precovering}.
	\end{equation}
	
	Dually, the adjunction $e\dashv r$ sends a monic $\mathcal{D}$-preenvelope of $eM$ to a monic $\mathcal{L}$-preenvelope of $M$.
	Conversely, since $s$ is additive and
	$sr\cong\mathrm{id}_{\mathcal{B}}$, we have
	\[
	r\mathcal{D}\subseteq\mathcal{L},
	\qquad
	s\mathcal{L}\subseteq
	\operatorname{add}(sr\mathcal{D})=\mathcal{D}.
	\]
	The dual assertion of Lemma~\ref{lem:transfer} therefore gives
	\begin{equation}\label{eq:preenvelop-transfer}
		\mathcal{D}\text{ is monic preenveloping}
		\quad\Longleftrightarrow\quad
		\mathcal{L}\text{ is monic preenveloping}.
	\end{equation}
	
	For every short exact sequence $\xi$ in $\mathcal{A}$, by the adjunctions $l\dashv e\dashv r$, we have
	\begin{equation}\label{eq:Hom-transfer}
		\begin{aligned}
			\operatorname{Hom}_{\mathcal{A}}(\mathcal{F},\xi)
			\text{ is exact}
			&\quad\Longleftrightarrow\quad
			\operatorname{Hom}_{\mathcal{B}}(\mathcal{C},e\xi)
			\text{ is exact},\\
			\operatorname{Hom}_{\mathcal{A}}(\xi,\mathcal{L})
			\text{ is exact}
			&\quad\Longleftrightarrow\quad
			\operatorname{Hom}_{\mathcal{B}}(e\xi,\mathcal{D})
			\text{ is exact}.
		\end{aligned}
	\end{equation}
	If $(\mathcal{C},\mathcal{D})$ satisfies
	\eqref{eq:balance-criterion}, then the exactness of $e$ allows us to apply this criterion to $e\xi$. By \eqref{eq:Hom-transfer}, $(\mathcal{F},\mathcal{L})$ satisfies the same criterion.
	
	Conversely, suppose that $(\mathcal{F},\mathcal{L})$ satisfies \eqref{eq:balance-criterion}. For any short exact sequence $\sigma$ in $\mathcal{B}$, the exactness of $i$ implies that $i\sigma$ is a short exact sequence. Since $ei\sigma\cong\sigma$, we obtain
	\[
	\begin{aligned}
		\operatorname{Hom}_{\mathcal{B}}(\mathcal{C},\sigma)
		\text{ is exact}
		&\Longleftrightarrow
		\operatorname{Hom}_{\mathcal{A}}(\mathcal{F},i\sigma)
		\text{ is exact}\\
		&\Longleftrightarrow
		\operatorname{Hom}_{\mathcal{A}}(i\sigma,\mathcal{L})
		\text{ is exact}\\
		&\Longleftrightarrow
		\operatorname{Hom}_{\mathcal{B}}(\sigma,\mathcal{D})
		\text{ is exact}.
	\end{aligned}
	\]
	Thus $(\mathcal{C},\mathcal{D})$ satisfies
	\eqref{eq:balance-criterion}.
	
	The result follows from Lemma~\ref{lem:balance-criterion},
	\eqref{eq:precover-transfer}, and
	\eqref{eq:preenvelop-transfer}.
\end{proof}

\begin{remark}\label{rem:exact-structure}
	Let $(\mathcal{C},\mathcal{D})$ be an admissible balanced pair in $\mathcal{B}$, and let $\mathcal{E}_{\mathcal{B}}$ denote
	its associated exact structure. By Theorem~\ref{thm:induced},
	the pair $(\operatorname{add}(l\mathcal{C}),
	\operatorname{add}(r\mathcal{D}))$
	is an admissible balanced pair in $\mathcal{A}$.
	Denote its associated exact structure by
	$\mathcal{E}_{\mathcal{A}}$.
	Then \eqref{eq:Hom-transfer} yields
	\[
	\mathcal{E}_{\mathcal{A}}
	=e^{-1}(\mathcal{E}_{\mathcal{B}})
	:=\{\xi\text{ short exact in }\mathcal{A}
	\mid e\xi\in\mathcal{E}_{\mathcal{B}}\}.
	\]
\end{remark}

\begin{example}\label{ex:section}
	The coinduction functor $r$ in Theorem~\ref{thm:induced} cannot in general be replaced by the section functor $i$. Let $k$ be a field and let
	$S=k[\varepsilon]/(\varepsilon^2)$.
	The canonical inclusion $k\to S$ and the augmentation
	$S\to k$, $\varepsilon\mapsto 0$, give a cleft extension with
	\[
	\mathcal{B}=k\text{-}\operatorname{Mod},
	\qquad
	\mathcal{A}=S\text{-}\operatorname{Mod},
	\]
	where $e$ is restriction of scalars along $k\to S$, $l=S\otimes_k-$, and $i$ is restriction of scalars along the augmentation $S\to k$. Moreover, $e$ has the right adjoint $r=\operatorname{Hom}_k(S,-)$.
	
	Since every short exact sequence of $k$-vector spaces splits, $(\mathcal{B},\mathcal{B})$ is an admissible balanced pair. However, every module in $i\mathcal{B}$ is annihilated by $\varepsilon$. For any $N\in i\mathcal{B}$ and any $S$-linear map $f:S\to N$, we have
	$f(\varepsilon)=\varepsilon f(1)=0$.
	Since $\varepsilon\neq 0$ in $S$, no such map is a monomorphism. Thus the regular left $S$-module admits no monomorphism into an object of $i\mathcal{B}$, so $i\mathcal{B}$ is not monic preenveloping. Consequently,
	$(\operatorname{add}(l\mathcal{B}),i\mathcal{B})$ is not an admissible balanced pair.
\end{example}

\begin{example}\label{ex:summands}
The essential image $l\mathcal{B}$ need not be closed under direct summands, even when
$(\mathcal{B}, \mathcal{B})$ is admissible balanced. Let $\mathcal{B} = k\text{-}\operatorname{Mod}$ for a field
$k$, and take $F = \mathrm{id}$ and $\eta = \mathrm{id}$. An object of
$\mathcal{A} = \mathcal{B} \ltimes_\eta F$ is a vector space with an idempotent endomorphism,
so that $\mathcal{A} \simeq k\text{-}\operatorname{Mod} \times k\text{-}\operatorname{Mod}$ and every object of
$\mathcal{A}$ is projective. By \eqref{eq:model},
\[
 lX = \left(X \oplus X, \begin{pmatrix}0 & 0\\ \mathrm{id} & \mathrm{id}\end{pmatrix}\right)
\cong (X, \mathrm{id}) \oplus (X, 0)
\]
for every $X \in \mathcal{B}$. Hence the direct summand $(k, \mathrm{id})$ of $lk$ does not
belong to $l\mathcal{B}$.
\end{example}

\section{Balanced pairs and cotorsion pairs under cleft extensions}\label{sec:structural}

In Section~\ref{sec:induced}, we compared admissible balance in $\mathcal{B}$ with
that of the induced and coinduced classes in $\mathcal{A}$. We now ask when balance
can instead be compared with that of the special classes $\Delta(\mathcal{C})$
and $\nabla(\mathcal{D})$.

We first give a criterion in terms of lifted complete hereditary cotorsion
pairs. We then show that this criterion applies to finite directed
extensions when $L_1F(\mathcal{C}) = 0$ and $R^1G(\mathcal{D}) = 0$, which leads to the
comparison in Theorem~C.

Throughout this section, let $(\mathcal{B},\mathcal{A},i,e,l)$ be a cleft extension of abelian categories. We assume that $e$
admits a right adjoint $r$ and that $\mathcal{B}$ has enough projectives and enough injectives.
\subsection{Special classes and cotorsion triplets}
For subcategories $\mathcal{C},\mathcal{D}\subseteq\mathcal{B}$,
put
\begin{equation}\label{eq:structural-classes}
	\begin{split}
		\Delta(\mathcal{C})
		&=\{M\in\mathcal{A}\mid L_1q(M)=0,\ qM\in\mathcal{C}\},\\
		\nabla(\mathcal{D})
		&=\{N\in\mathcal{A}\mid R^1s(N)=0,\ sN\in\mathcal{D}\}.
	\end{split}
\end{equation}
By \eqref{eq:L1q-model} and \eqref{eq:R1s-model},
these agree, in the corresponding models, with the
classes defined in \cite[Subsections~3.1 and~3.6]{Hu}. We have
\begin{equation}\label{eq:induced-in-structural}
	l\mathcal{C}\subseteq\Delta(\mathcal{C}),
	\qquad
	r\mathcal{D}\subseteq\nabla(\mathcal{D}).
\end{equation}
Indeed, let $X\in\mathcal{C}$ and let $I$ be an injective object of $\mathcal{B}$. Consider a short exact sequence
\[
0 \to iI \to E\xrightarrow{p}lX \to 0.
\]
Since $e$ is exact and $eiI\cong I$ is injective,
its image under $e$ splits. Choose a section
$t:elX\to eE$ of $e(p)$, and put
$h=t\circ\nu_X:X\to eE$, where $\nu$ is the unit
of $l\dashv e$. Then $e(p)\circ h=\nu_X$.
By naturality of the adjunction isomorphism, the adjoint
$\widetilde h = \mu_E\circ l(h):lX\to E$ satisfies
$p\circ\widetilde h=\mathrm{id}_{lX}$, where $\mu$ is the counit of $l\dashv e$. Consequently,
$\operatorname{Ext}_{\mathcal{A}}^1(lX,iI)=0$.
By \eqref{eq:injective-test},
\[
\operatorname{Hom}_{\mathcal{B}}(L_1q(lX),I)=0
\]
for every injective $I\in\mathcal{B}$.
Since $\mathcal{B}$ has enough injective objects,
$L_1q(lX)=0$. Together with $qlX\cong X$, this proves $lX\in\Delta(\mathcal{C})$. Dually, we have $r\mathcal{D}\subseteq\nabla(\mathcal{D})$.

\begin{lemma}\label{lem:reflection}
	Let $\mathcal{C},\mathcal{D}\subseteq\mathcal{B}$ and
	$\mathcal{F},\mathcal{L}\subseteq\mathcal{A}$ be subcategories such that
	\[
	l\mathcal{C}\subseteq\mathcal{F},\quad
	q\mathcal{F}\subseteq\mathcal{C},\qquad
	r\mathcal{D}\subseteq\mathcal{L},\quad
	s\mathcal{L}\subseteq\mathcal{D}.
	\]
	If $(\mathcal{F},\mathcal{L})$ is an admissible balanced pair in $\mathcal{A}$, then $(\mathcal{C},\mathcal{D})$ is an
	admissible balanced pair in $\mathcal{B}$.
	In particular, if $(\Delta(\mathcal{C}),\nabla(\mathcal{D}))$
	is an admissible balanced pair in $\mathcal{A}$, then $(\mathcal{C},\mathcal{D})$ is an admissible balanced pair in $\mathcal{B}$.
\end{lemma}

\begin{proof}
	For each $B\in\mathcal{B}$, choose an epic
	$\mathcal{F}$-precover and a monic $\mathcal{L}$-preenvelope of $iB$. By Lemma~\ref{lem:transfer}, these induce an epic
	$\mathcal{C}$-precover and a monic $\mathcal{D}$-preenvelope of $B$, respectively. Thus $\mathcal{C}$ is epic precovering and $\mathcal{D}$ is monic preenveloping.
	
	Since $ql\cong\mathrm{id}_{\mathcal{B}}$ and
	$sr\cong\mathrm{id}_{\mathcal{B}}$, the inclusions in the hypotheses give $q\mathcal{F}=\mathcal{C}$ and
	$s\mathcal{L}=\mathcal{D}$. Let $\sigma$ be a short exact sequence in $\mathcal{B}$.
	Since $i$ is exact, $i\sigma$ is short exact in $\mathcal{A}$. Using $q\dashv i$, the balance condition for
	$(\mathcal{F},\mathcal{L})$, and $i\dashv s$, we obtain
	\[
	\begin{aligned}
		\operatorname{Hom}_{\mathcal{B}}(\mathcal{C},\sigma)
		\text{ is exact}
		&\Longleftrightarrow
		\operatorname{Hom}_{\mathcal{A}}(\mathcal{F},i\sigma)
		\text{ is exact}\\
		&\Longleftrightarrow
		\operatorname{Hom}_{\mathcal{A}}(i\sigma,\mathcal{L})
		\text{ is exact}\\
		&\Longleftrightarrow
		\operatorname{Hom}_{\mathcal{B}}(\sigma,\mathcal{D})
		\text{ is exact}.
	\end{aligned}
	\]
	Hence $(\mathcal{C},\mathcal{D})$ is an admissible balanced pair by Lemma~\ref{lem:balance-criterion}.
	
	For the final assertion, \eqref{eq:induced-in-structural} gives $l\mathcal{C}\subseteq\Delta(\mathcal{C})$ and
	$r\mathcal{D}\subseteq\nabla(\mathcal{D})$, while
	$q\Delta(\mathcal{C})\subseteq\mathcal{C}$ and
	$s\nabla(\mathcal{D})\subseteq\mathcal{D}$ follow from \eqref{eq:structural-classes}.
\end{proof}

Lemma~\ref{lem:reflection} shows that if
$(\Delta(\mathcal{C}),\nabla(\mathcal{D}))$ is an admissible balanced
pair in $\mathcal{A}$, then $(\mathcal{C},\mathcal{D})$ is an
admissible balanced pair in $\mathcal{B}$. The next theorem gives
conditions under which the converse also holds.

\begin{theorem}\label{thm:structural-criterion}
Let $(\mathcal{C}, \mathcal{Y})$ and $(\mathcal{X}, \mathcal{D})$ be complete hereditary cotorsion pairs in $\mathcal{B}$
such that
\begin{equation}\label{eq:core-compatibility}
 \mathcal{C} \cap \mathcal{Y} \subseteq \mathcal{X}, \qquad \mathcal{X} \cap \mathcal{D} \subseteq \mathcal{Y}.
\end{equation}
Suppose that
\[
 (\Delta(\mathcal{C}), e^{-1}(\mathcal{Y})), \qquad
 (e^{-1}(\mathcal{X}), \nabla(\mathcal{D}))
\]
are complete hereditary cotorsion pairs in $\mathcal{A}$. Then the following
conditions are equivalent:
\begin{enumerate}
 \item $(\mathcal{C}, \mathcal{D})$ is admissible balanced;
 \item $(\Delta(\mathcal{C}), \nabla(\mathcal{D}))$ is admissible balanced;
 \item $\mathcal{Y} = \mathcal{X}$;
 \item $\Delta(\mathcal{C})^\perp = {}^\perp\nabla(\mathcal{D})$.
\end{enumerate}
\end{theorem}

\begin{proof}
	Suppose that (1) holds. Then $(\mathcal{C},\mathcal{D})$
	is an admissible balanced pair in $\mathcal{B}$. Since
	$(\mathcal{C},\mathcal{Y})$ and $(\mathcal{X},\mathcal{D})$
	are complete hereditary cotorsion pairs satisfying
	\eqref{eq:core-compatibility}, \cite[Corollary~4.8]{EPZ}
	implies that $\mathcal{Y}=\mathcal{X}$. Thus $(1)\Rightarrow(3)$.
	
	Suppose next that (3) holds. Then $e^{-1}(\mathcal{Y})=e^{-1}(\mathcal{X})$. Therefore,
	$(\Delta(\mathcal{C}),e^{-1}(\mathcal{Y}),
	\nabla(\mathcal{D}))$ is a complete hereditary cotorsion triplet in $\mathcal{A}$.
	By Lemma~\ref{lem:adjoints}, $\mathcal{A}$ has enough
	projectives and enough injectives. Hence
	Lemma~\ref{lem:triplet-balanced} gives (2).
	The implication $(2)\Rightarrow(1)$ follows from
	Lemma~\ref{lem:reflection}.
	
	Finally, since $(\Delta(\mathcal{C}), e^{-1}(\mathcal{Y}))$ and
	$(e^{-1}(\mathcal{X}), \nabla(\mathcal{D}))$ are complete hereditary cotorsion pairs, we have
	\[
	\Delta(\mathcal{C})^\perp=e^{-1}(\mathcal{Y}),
	\qquad
	{}^\perp\nabla(\mathcal{D})=e^{-1}(\mathcal{X}).
	\]
	Since $ei\cong\mathrm{id}_{\mathcal{B}}$, for every $B\in\mathcal{B}$ we have
	\[
	B\in\mathcal{Y}\iff iB\in e^{-1}(\mathcal{Y}),
	\qquad
	B\in\mathcal{X}\iff iB\in e^{-1}(\mathcal{X}).
	\]
	Thus $e^{-1}(\mathcal{Y})=e^{-1}(\mathcal{X})$ if and only if $\mathcal{Y}=\mathcal{X}$. This proves $(3)\Leftrightarrow(4)$ and completes the proof.
\end{proof}

The first equality below is due to Hu \cite[Lemma~3.8]{Hu} when $F^2=0$ and to Chen, Qin and Yang \cite[Lemma~3.3]{ChenQinYang} for general $\eta$-extensions, and we include a short proof using \eqref{eq:five-term} and \eqref{eq:injective-test}.

\begin{lemma}\label{lem:structural-orthogonality}
	If $\mathcal{Y}$ contains all injective objects of $\mathcal{B}$, then
	$\Delta({}^\perp\mathcal{Y}) = {}^\perp(i\mathcal{Y})$.
	Dually, if $\mathcal{X}$ contains all projective objects of $\mathcal{B}$, then
	$\nabla(\mathcal{X}^\perp) = (i\mathcal{X})^\perp$.
\end{lemma}

\begin{proof}
	Let $M\in\Delta({}^\perp\mathcal{Y})$. Then $L_1q(M)=0$ and $qM\in{}^\perp\mathcal{Y}$. For every $Y\in\mathcal{Y}$,
	\eqref{eq:five-term} gives
	\[
	\operatorname{Ext}^1_{\mathcal{A}}(M,iY)
	\cong \operatorname{Ext}^1_{\mathcal{B}}(qM,Y)=0.
	\]
	Hence $M\in{}^\perp(i\mathcal{Y})$.
	
	Conversely, let $M\in{}^\perp(i\mathcal{Y})$. For every $Y\in\mathcal{Y}$, the first map in \eqref{eq:five-term} gives an
	injection
	\[
	\operatorname{Ext}^1_{\mathcal{B}}(qM,Y)
	\hookrightarrow \operatorname{Ext}^1_{\mathcal{A}}(M,iY)=0.
	\]
	Thus $qM\in{}^\perp\mathcal{Y}$. Moreover, every injective object $I$ of $\mathcal{B}$ belongs to $\mathcal{Y}$, so
	\eqref{eq:injective-test} gives
	\[
	\operatorname{Hom}_{\mathcal{B}}(L_1q(M),I)
	\cong \operatorname{Ext}^1_{\mathcal{A}}(M,iI)=0.
	\]
	Since $\mathcal{B}$ has enough injectives, $L_1q(M)$ embeds into an injective object. Hence $L_1q(M)=0$. This proves
	$M\in\Delta({}^\perp\mathcal{Y})$. The second equality follows by duality.
\end{proof}

\begin{corollary}\label{cor:any-eta}
Let $(\mathcal{C}, \mathcal{Y})$ be a cotorsion pair in $\mathcal{B}$.
\begin{enumerate}
 \item $\Delta(\mathcal{C}) = {}^\perp(i\mathcal{Y})$, and $(\Delta(\mathcal{C}), \Delta(\mathcal{C})^\perp)$
 is a cotorsion pair in $\mathcal{A}$.
 \item ${}^\perp e^{-1}(\mathcal{Y}) \subseteq \Delta(\mathcal{C})$, with equality whenever
 every object of $e^{-1}(\mathcal{Y})$ has a finite filtration with factors in
 $i\mathcal{Y}$.
 \item  $({}^\perp e^{-1}(\mathcal{Y}), e^{-1}(\mathcal{Y}))$ is a
 cotorsion pair in $\mathcal{A}$ if and only if $L_1F(\mathcal{C}) = 0$. In this case
 $\Delta(\mathcal{C})^\perp \subseteq e^{-1}(\mathcal{Y})$.
\end{enumerate}
\end{corollary}

\begin{proof}
	Since $(\mathcal{C},\mathcal{Y})$ is a cotorsion pair, $\mathcal{Y}$ contains all injective objects of $\mathcal{B}$ and
	$\mathcal{C}={}^\perp\mathcal{Y}$.
	Part~(1) follows from Lemma~\ref{lem:structural-orthogonality}; see also
	\cite[Lemma~3.3 and Theorem~3.4]{ChenQinYang}.
	
	Now we prove part (2). Since $ei\cong\mathrm{id}_{\mathcal{B}}$, we have
	$i\mathcal{Y}\subseteq e^{-1}(\mathcal{Y})$. Consequently,
	\[
	{}^\perp e^{-1}(\mathcal{Y})
	\subseteq {}^\perp(i\mathcal{Y})
	=\Delta(\mathcal{C}).
	\]
	Suppose now that the filtration condition in (2) holds.
	Let $M\in\Delta(\mathcal{C})$ and $N\in e^{-1}(\mathcal{Y})$,
	and choose a finite filtration
	\[
	0=N_0\subseteq N_1\subseteq\cdots\subseteq N_m=N
	\]
	with $N_j/N_{j-1}\in i\mathcal{Y}$ for every $j$.
	By (1),
	$\operatorname{Ext}^1_{\mathcal{A}}(M,N_j/N_{j-1})=0$.
	Applying $\operatorname{Hom}_{\mathcal{A}}(M,-)$ to
	\[
	0\longrightarrow N_{j-1}\longrightarrow N_j
	\longrightarrow N_j/N_{j-1}\longrightarrow 0
	\]
	 gives an exact sequence
	\[
	\operatorname{Ext}^1_{\mathcal{A}}(M,N_{j-1})
	\longrightarrow \operatorname{Ext}^1_{\mathcal{A}}(M,N_j)
	\longrightarrow \operatorname{Ext}^1_{\mathcal{A}}(M,N_j/N_{j-1}).
	\]
	Since the last term vanishes and $N_0=0$, induction gives
	$\operatorname{Ext}^1_{\mathcal{A}}(M,N_j)=0$ for every $j$. Hence $\operatorname{Ext}^1_{\mathcal{A}}(M,N)=0$, proving $\Delta(\mathcal{C})\subseteq{}^\perp e^{-1}(\mathcal{Y})$.
	
	For (3), if $L_1F(\mathcal{C}) = 0$, then $({}^\perp e^{-1}(\mathcal{Y}), e^{-1}(\mathcal{Y}))$ is a cotorsion pair in $\mathcal{A}$ by \cite[Theorem~3.10]{ChenQinYang}; see \cite[Theorem~3.5]{Hu} for nilpotent $\eta$. Moreover, the inclusion $\Delta(\mathcal{C})^\perp \subseteq e^{-1}(\mathcal{Y})$ follows from \cite[Proposition~3.11]{ChenQinYang}.
	
	Conversely, suppose that $\bigl({}^\perp e^{-1}(\mathcal{Y}),e^{-1}(\mathcal{Y})\bigr)$
	is a cotorsion pair in $\mathcal{A}$. Let $I$ be an injective object of $\mathcal{B}$.
	Since $e$ is exact and $e\dashv r$, the object $rI$ is injective in $\mathcal{A}$. Hence $rI\in e^{-1}(\mathcal{Y})$, and therefore $I\oplus GI\cong erI\in\mathcal{Y}$. As $\mathcal{Y}$ is closed under direct summands, $GI\in\mathcal{Y}$.
	
	For $C\in\mathcal{C}$, choose a projective resolution $P_\bullet\to C$. The adjunction $F\dashv G$ and the injectivity of $I$ yield
	\[
	0 = \operatorname{Ext}^1_{\mathcal{B}}(C, GI) \cong H^1\operatorname{Hom}_{\mathcal{B}}(FP_\bullet, I)
	\cong \operatorname{Hom}_{\mathcal{B}}(L_1F(C), I),
	\]
	where the last isomorphism follows from the exactness of $\operatorname{Hom}_{\mathcal{B}}(-,I)$.
	Since this holds for every injective object $I$ and $\mathcal{B}$ has enough injective objects, $L_1F(C)=0$. Thus $L_1F(\mathcal{C})=0$.
\end{proof}

\subsection{Finite directed extensions}\label{subsec:directed}

Let $\mathcal{B}=\prod_{j=1}^n\mathcal{B}_j$, with coordinate embeddings $\iota_j:\mathcal{B}_j\to\mathcal{B}$ and projections
$\pi_j:\mathcal{B}\to\mathcal{B}_j$. Set
\[
F_{ij}=\pi_iF\iota_j,
\qquad
G_{ji}=\pi_jG\iota_i.
\]
Since $F\dashv G$ and $\iota_j\dashv\pi_j\dashv\iota_j$,
we have $F_{ij}\dashv G_{ji}$. Moreover, for
$X=(X_1,\ldots,X_n)\in\mathcal{B}$, there are natural isomorphisms
\[
(FX)_i\cong\bigoplus_{j=1}^nF_{ij}X_j,
\qquad
(GX)_j\cong\bigoplus_{i=1}^nG_{ji}X_i.
\]

\begin{definition}\label{def:directed}
	The cleft extension $\mathcal{A}=\mathcal{B}\ltimes_\eta F$
	is called \emph{finite directed} with respect to the given
	product decomposition if
	\begin{equation}\label{eq:directed}
		F_{ij}=0\qquad\text{whenever }i\leq j.
	\end{equation}
\end{definition}

For such an extension, we have
\[
(FX)_i\cong\bigoplus_{j<i}F_{ij}X_j,
\qquad
(GX)_j\cong\bigoplus_{i>j}G_{ji}X_i.
\]
Iterating these formulas gives $F^n=0$ and $G^n=0$.

The aim of this subsection is to lift cotorsion pairs from the categories
$\mathcal{B}_j$ to $\mathcal{A}$. We fix the following notation. The extension
$\mathcal{A}=\mathcal{B}\ltimes_\eta F$ is finite directed, with $F\dashv G$ and
$\eta:F^2\to F$ any associative natural transformation, and each
$\mathcal{B}_j$ has enough projective and injective objects. For each $j$, let
$(\mathcal{C}_j,\mathcal{Y}_j)$ be a cotorsion pair in $\mathcal{B}_j$. The product classes
$\mathcal{C}$ and $\mathcal{Y}$ form a cotorsion pair $(\mathcal{C},\mathcal{Y})$ in $\mathcal{B}$, and we
put
\[
\mathcal{S}=\Delta(\mathcal{C}),\qquad \mathcal{T}=e^{-1}(\mathcal{Y}).
\]

The first step is a filtration of the objects of $\mathcal{A}$ by the objects
supported at a single coordinate.

\begin{lemma}\label{lem:vertex-filtration}
	Every $M=((M_a)_{a=1}^n,f)\in\mathcal{A}$ has a filtration
	\begin{equation}\label{eq:vertex-filtration}
		0=M_{\geq n+1}\subseteq M_{\geq n}\subseteq\cdots\subseteq M_{\geq1}=M,
		\qquad M_{\geq j}/M_{\geq j+1}\cong i\iota_j(M_j),
	\end{equation}
	where $M_{\geq j}$ has underlying object $(0,\ldots,0,M_j,\ldots,M_n)$ and
	the structure map induced by $f$.
\end{lemma}

\begin{proof}
	Let $X\in\mathcal{B}$ with $X_k=0$ for all $k<j$. By the directed condition,
	$(FX)_k=0$ for all $k\leq j$. Hence $f$ restricts to a structure map on
	$(0,\ldots,0,M_j,\ldots,M_n)$, so that $M_{\geq j}$ is a subobject of $M$.
	The quotient $M_{\geq j}/M_{\geq j+1}$ has underlying object $\iota_j(M_j)$.
	Its structure map is zero because $F\iota_j(M_j)$ vanishes at the $j$th
	coordinate. Thus $M_{\geq j}/M_{\geq j+1}\cong i\iota_j(M_j)$.
\end{proof}

Using this filtration, we prove that $\mathcal{S}={}^\perp\mathcal{T}$, and determine when
$(\mathcal{S},\mathcal{T})$ is a hereditary cotorsion pair.

\begin{proposition}\label{prop:directed-pair}
	\begin{enumerate}
		\item $\mathcal{S}={}^\perp\mathcal{T}$.
		\item $(\mathcal{S},\mathcal{T})$ is a cotorsion pair in $\mathcal{A}$ if and only if
		$L_1F(\mathcal{C})=0$.
		\item If $L_1F(\mathcal{C})=0$ and every $(\mathcal{C}_j,\mathcal{Y}_j)$ is hereditary, then
		$(\mathcal{S},\mathcal{T})$ is hereditary.
	\end{enumerate}
\end{proposition}

\begin{proof}
	(1) Let $M\in\mathcal{T}$. Each $M_j$ belongs to $\mathcal{Y}_j$, so that the factors
	$i\iota_j(M_j)$ of the filtration \eqref{eq:vertex-filtration} belong to
	$i\mathcal{Y}$. The assertion follows from Corollary~\ref{cor:any-eta}(2).

	(2) By (1), the pair $(\mathcal{S},\mathcal{T})$ is
	$({}^\perp e^{-1}(\mathcal{Y}),e^{-1}(\mathcal{Y}))$, and the assertion is
	Corollary~\ref{cor:any-eta}(3).

	(3) Since each $(\mathcal{C}_j,\mathcal{Y}_j)$ is hereditary, the class $\mathcal{Y}$ is closed
	under cokernels of monomorphisms. Hence $\mathcal{T}$ is also closed under
	cokernels of monomorphisms, since $e$ is exact. Moreover, since
	$(\mathcal{S},\mathcal{T})$ is a cotorsion pair by (2), $\mathcal{T}$ contains all injective
	objects of $\mathcal{A}$, and $\mathcal{A}$ has enough injective objects by
	Lemma~\ref{lem:adjoints}. Thus dimension shifting yields
	$\operatorname{Ext}^m_{\mathcal{A}}(\mathcal{S},\mathcal{T})=0$ for all $m\geq1$, so that
	$(\mathcal{S},\mathcal{T})$ is hereditary.
\end{proof}

We now turn to the completeness of $(\mathcal{S},\mathcal{T})$, which will be proved
by induction on $n$. The induction uses the extensions obtained by
discarding the first coordinates, which we introduce first.

Fix $1\leq j<n$ and put
$\mathcal{B}_{>j}=\prod_{k>j}\mathcal{B}_k$.
Let $\iota^{\mathcal{B}}_{>j}$ and $\pi_{>j}$ denote the coordinate embedding and projection, respectively. By the directed condition, for any object
$X=(X_1,\ldots,X_n)\in\mathcal{B}$ with $X_k=0$ for $k\leq j$, we also have $(FX)_k=0$ for $k\leq j$. Hence $F$ and $\eta$ restrict to $\mathcal{B}_{>j}$, giving an extension
\[
\mathcal{A}_{>j}
=\mathcal{B}_{>j}\ltimes_{\eta_{>j}}F_{>j},
\qquad
F_{>j}=\pi_{>j}F\iota^{\mathcal{B}}_{>j},
\]
where $\eta_{>j}$ is induced by $\eta$.
The functor $F_{>j}$ has right adjoint
$G_{>j}=\pi_{>j}G\iota^{\mathcal{B}}_{>j}$.
This extension is finite directed.

Adding zero objects in the first $j$ coordinates defines an exact functor
$\iota^{\mathcal{A}}_{>j}:\mathcal{A}_{>j}\to\mathcal{A}$.
Since adding zero coordinates commutes with taking kernels, images and cokernels, the formulas for $q$ and $L_1q$ give
\[
q\iota^{\mathcal{A}}_{>j}
\cong\iota^{\mathcal{B}}_{>j}q_{>j},
\qquad
L_1q\,\iota^{\mathcal{A}}_{>j}
\cong\iota^{\mathcal{B}}_{>j}L_1q_{>j}.
\]
Let
$N\in\Delta_{>j}(\prod_{k>j}\mathcal{C}_k)$.
By definition, $L_1q_{>j}(N)=0$ and
$q_{>j}(N)\in\prod_{k>j}\mathcal{C}_k$.
The above isomorphisms give
\[
L_1q(\iota^{\mathcal{A}}_{>j}N)
\cong\iota^{\mathcal{B}}_{>j}L_1q_{>j}(N)=0
\]
and
\[
q(\iota^{\mathcal{A}}_{>j}N)
\cong\iota^{\mathcal{B}}_{>j}q_{>j}(N).
\]
The first $j$ coordinates of
$\iota^{\mathcal{B}}_{>j}q_{>j}(N)$ are zero, and its $k$th coordinate belongs to $\mathcal{C}_k$ for every $k>j$. Since $0\in\mathcal{C}_k$ for $k\leq j$, we have
$\iota^{\mathcal{B}}_{>j}q_{>j}(N)\in\mathcal{C}$. Thus $\iota^{\mathcal{A}}_{>j}N\in
\Delta(\mathcal{C})=\mathcal{S}$, which means that $\iota^{\mathcal{A}}_{>j}
\Bigl(\Delta_{>j}\bigl(\prod_{k>j}\mathcal{C}_k\bigr)\Bigr) \subseteq\mathcal{S}$. Let  $e_{>j}$ denote the forgetful functor for $\mathcal{A}_{>j}$. Likewise, let
$N\in e_{>j}^{-1}(\prod_{k>j}\mathcal{Y}_k)$.
Then $e(\iota^{\mathcal{A}}_{>j}N)
\cong\iota^{\mathcal{B}}_{>j}e_{>j}(N)$.
The first $j$ coordinates of this object are zero, and its $k$th coordinate belongs to $\mathcal{Y}_k$ for every $k>j$. Since each $\mathcal{Y}_k$ contains the zero object,
$e(\iota^{\mathcal{A}}_{>j}N)\in\mathcal{Y}$.
Thus $\iota^{\mathcal{A}}_{>j}N\in
e^{-1}(\mathcal{Y})=\mathcal{T}$, which means that $\iota^{\mathcal{A}}_{>j} \Bigl(e_{>j}^{-1}\bigl(\prod_{k>j}\mathcal{Y}_k\bigr)\Bigr) \subseteq\mathcal{T}$.

The next lemma is the key step of the induction. It provides special
$\mathcal{S}$-precovers for the objects $i\iota_j(X)$ with $X\in\mathcal{B}_j$.

\begin{lemma}\label{lem:vertex-precover}
	Suppose that $L_1F(\mathcal{C})=0$, so that $(\mathcal{S},\mathcal{T})$ is a cotorsion pair
	by Proposition~\ref{prop:directed-pair}(2). Let $1\leq j\leq n$ be such
	that $(\mathcal{C}_j,\mathcal{Y}_j)$ is complete and, if $j<n$, such that
	\[
	\Bigl(\Delta_{>j}\bigl(\prod_{k>j}\mathcal{C}_k\bigr),\
	e_{>j}^{-1}\bigl(\prod_{k>j}\mathcal{Y}_k\bigr)\Bigr)
	\]
	is a complete cotorsion pair in $\mathcal{A}_{>j}$. Then $i\iota_jX$ has a
	special $\mathcal{S}$-precover for every $X\in\mathcal{B}_j$.
\end{lemma}

\begin{proof}
	For $X\in\mathcal{B}_j$, the completeness of
	$(\mathcal{C}_j,\mathcal{Y}_j)$ gives a special
	$\mathcal{C}_j$-precover sequence
	\[
	0\longrightarrow Y\longrightarrow C
	\xrightarrow{\;\;a\;\;}X\longrightarrow0,
	\qquad C\in\mathcal{C}_j,\quad Y\in\mathcal{Y}_j.
	\]
	Let $\widetilde a:l\iota_jC\to i\iota_jX$ be the morphism corresponding to $\iota_ja$ under $l\dashv e$. Under the decomposition
	$el\iota_jC\cong\iota_jC\oplus F\iota_jC$, we have $e(\widetilde a)=(\iota_ja,0): \iota_jC \oplus F\iota_jC \longrightarrow \iota_jX$.
	Indeed, write $e(\widetilde a)=(u,v)$.
	By adjunction, $e(\widetilde a)\nu_{\iota_jC}=\iota_ja$. Since $\nu_{\iota_jC}=\binom{\mathrm{id}}{0}$ under the above
	decomposition, we obtain $u=\iota_ja$.
	Moreover,
	\[
	l\iota_jC=
	\left(
	\iota_jC\oplus F\iota_jC,\,
	\begin{pmatrix}
		0&0\\
		\mathrm{id}_{F\iota_jC}&\eta_{\iota_jC}
	\end{pmatrix}
	\right),
	\qquad
	i\iota_jX=(\iota_jX,0).
	\]
	As $\widetilde a$ is a morphism in $\mathcal{A}$,
	compatibility with these structure maps gives
	\[
	(u,v)
	\begin{pmatrix}
		0&0\\
		\mathrm{id}&\eta_{\iota_jC}
	\end{pmatrix}
	=(v,v\eta_{\iota_jC})=0.
	\]
	Thus $v=0$, and hence $e(\widetilde a)=(\iota_ja,0)$.

	Since $a$ is an epimorphism, so is
	$e(\widetilde a)=(\iota_ja,0)$.
	As $e$ reflects epimorphisms, $\widetilde a$ is an epimorphism. Put $K=\operatorname{Ker}\widetilde a$. Since $e$ is exact, $eK\cong\operatorname{Ker}(\iota_ja,0)
	\cong\iota_jY\oplus F\iota_jC$.
	By the directed condition, $(F\iota_jC)_k=0$ for $k\leq j$.
	Hence
	\[
	(eK)_j\cong Y,\qquad (eK)_k=0\quad(k<j).
	\]

	If $j=n$, the directed condition gives $F\iota_nC=0$.
	Hence $l\iota_nC\cong i\iota_nC$ and
	$K\cong i\iota_nY\in\mathcal{T}$.

	Since $\iota_nC\in\mathcal{C}$ and
	$l\mathcal{C}\subseteq\Delta(\mathcal{C})=\mathcal{S}$, we have $l\iota_nC\in\mathcal{S}$.
	Then the exact sequence
	\[
	0\longrightarrow K\longrightarrow l\iota_nC
	\xrightarrow{\widetilde a}i\iota_nX\longrightarrow0
	\]
	is a special $\mathcal{S}$-precover sequence.

	We may therefore assume $j<n$. Since $(eK)_k=0$ for $k<j$ and $(eK)_j\cong Y$,
	the coordinate filtration gives an exact sequence
	\[
	0\longrightarrow K_{>j}\longrightarrow K
	\longrightarrow i\iota_jY\longrightarrow0,
	\]
	where $K_{>j}=K_{\geq j+1}$ is viewed as an object of $\mathcal{A}_{>j}$ via the embedding
	$\iota^{\mathcal{A}}_{>j}$. By hypothesis, choose in $\mathcal{A}_{>j}$ a special preenvelope sequence
	\[
	0 \to K_{>j}\xrightarrow u T_{>j} \to S_{>j} \to 0,
	\qquad T_{>j} \in \mathcal{T}, \quad S_{>j} \in \mathcal{S}.
	\]
	Form the following two pushouts successively:
	\[
	\begin{tikzcd}[column sep=small]
		0\arrow[r]&K_{>j}\arrow[r]\arrow[d,"u"']&
		K\arrow[r]\arrow[d]&i\iota_jY\arrow[r]\arrow[d,equal]&0\\
		0\arrow[r]&T_{>j}\arrow[r]&K'\arrow[r]&i\iota_jY\arrow[r]&0
	\end{tikzcd}
	\qquad
	\begin{tikzcd}[column sep=small]
		0\arrow[r]&K\arrow[r]\arrow[d]&l\iota_jC\arrow[r]\arrow[d]&
		i\iota_jX\arrow[r]\arrow[d,equal]&0\\
		0\arrow[r]&K'\arrow[r]&P\arrow[r]&i\iota_jX\arrow[r]&0
	\end{tikzcd}
	\]
	The bottom row of the first diagram is exact.
	Since $T_{>j},i\iota_jY\in\mathcal{T}$ and
	$\mathcal{T}$ is closed under extensions, we have $K'\in\mathcal{T}$. Moreover, the induced map $K\to K'$ is a monomorphism with cokernel
	\[
	K'/K\cong\operatorname{Coker}u\cong S_{>j}.
	\]
	The second pushout therefore gives an exact sequence
	\[
	0\longrightarrow l\iota_jC\longrightarrow P
	\longrightarrow S_{>j}\longrightarrow0.
	\]
	Since $l\iota_jC,S_{>j}\in\mathcal{S}$ and
	$\mathcal{S}$ is closed under extensions, $P\in\mathcal{S}$. Thus the bottom row of the second diagram,
	\[
	0\longrightarrow K'\longrightarrow P
	\longrightarrow i\iota_jX\longrightarrow0,
	\]
	is the required special $\mathcal{S}$-precover sequence.
\end{proof}

For finite directed extensions, we establish a necessary and sufficient condition for the induced cotorsion pairs to be complete and hereditary.
\begin{theorem}\label{thm:directed}
	Let $\mathcal{A} = \mathcal{B} \ltimes_\eta F$ be finite directed, with $F \dashv G$ and
	$\eta: F^2 \to F$ any associative natural transformation, and suppose that
	each $\mathcal{B}_j$ has enough projective and injective objects. For each $j$,
	let $(\mathcal{C}_j, \mathcal{Y}_j)$ and $(\mathcal{X}_j, \mathcal{D}_j)$ be complete hereditary cotorsion
	pairs in $\mathcal{B}_j$, and denote the product classes by $\mathcal{C}, \mathcal{Y}, \mathcal{X}, \mathcal{D}$.
	\begin{enumerate}
		\item $(\Delta(\mathcal{C}), e^{-1}(\mathcal{Y}))$ is a complete hereditary cotorsion
		pair in $\mathcal{A}$ if and only if $L_1F(\mathcal{C}) = 0$.
		\item $(e^{-1}(\mathcal{X}), \nabla(\mathcal{D}))$ is a complete hereditary cotorsion
		pair in $\mathcal{A}$ if and only if $R^1G(\mathcal{D}) = 0$.
	\end{enumerate}
	The two vanishing conditions can be checked componentwise: they are
	respectively equivalent to $L_1F_{ij}(\mathcal{C}_j) = 0$ and
	$R^1G_{ji}(\mathcal{D}_i) = 0$ for all $i>j$.
\end{theorem}

\begin{proof}
	It suffices to prove (1). If $(\mathcal{S}, \mathcal{T})$ is a cotorsion pair, then $L_1F(\mathcal{C}) = 0$ by
	Corollary~\ref{cor:any-eta}(3). Conversely, suppose that $L_1F(\mathcal{C}) = 0$.
	Then $(\mathcal{S},\mathcal{T})$ is a hereditary cotorsion pair by
	Proposition~\ref{prop:directed-pair}.

	We prove completeness by induction on $n$.
	If $n=1$, the directed condition gives $F=0$.
	Thus $\mathcal{A}\simeq\mathcal{B}_1$, under which $(\mathcal{S},\mathcal{T})=(\mathcal{C}_1,\mathcal{Y}_1)$.
	Let $n>1$ and $1\leq j<n$.
	Coordinatewise projective resolutions show that
	$L_1F_{>j}\Bigl(\prod_{k>j}\mathcal{C}_k\Bigr)=0$. Thus the induction hypothesis applies to $\mathcal{A}_{>j}$.
	Hence Lemma~\ref{lem:vertex-precover} applies for every
	$1\leq j\leq n$.

	Thus every object $i\iota_jX$ with $X \in \mathcal{B}_j$ has a special
	$\mathcal{S}$-precover. Since $(\mathcal{S}, \mathcal{T})$ is a hereditary
	cotorsion pair, the class of objects of $\mathcal{A}$ admitting a special
	$\mathcal{S}$-precover is closed under extensions by \cite[Proposition~2.6(1)]{Hu}.
	For any $M\in\mathcal{A}$, the filtration
	\eqref{eq:vertex-filtration} gives exact sequences
	\[
	0\longrightarrow M_{\geq j+1}\longrightarrow M_{\geq j}
	\longrightarrow i\iota_j(M_j)\longrightarrow0.
	\]
	Starting with $M_{\geq n+1}=0$ and proceeding by descending induction on $j$, we conclude that $M=M_{\geq1}$ admits a special $\mathcal{S}$-precover. Since $\mathcal{A}$ has enough injective objects by Lemma~\ref{lem:adjoints} and every object admits a special $\mathcal{S}$-precover, the pullback argument in the proof of \cite[Proposition~7.1.7]{EnochsJenda}
	gives a special $\mathcal{T}$-preenvelope for every object of $\mathcal{A}$. Hence $(\mathcal{S},\mathcal{T})$ is complete.

	Statement (2) follows dually. The componentwise equivalences follow from
	coordinatewise resolutions.
\end{proof}

\begin{remark}\label{rem:hu-comparison}
For general $\eta$-extensions, the known completeness results for
$({}^\perp e^{-1}(\mathcal{Y}), e^{-1}(\mathcal{Y}))$ and for
$(\Delta(\mathcal{C}), \Delta(\mathcal{C})^\perp)$ require further conditions
such as $F(\mathcal{C}) \subseteq \mathcal{Y}$; see \cite[Theorem~3.13]{Hu},
\cite[Theorem~4.13]{LiuGao} and
\cite[Theorems~3.21 and~3.23]{ChenQinYang}. No such condition is needed in
Theorem~\ref{thm:directed}. For $n = 2$ the extension is a comma category,
and the sufficiency in parts (1) and (2) of Theorem~\ref{thm:directed} is proved in
\cite[Lemma~3.3 and Proposition~3.4]{HuZhuComma} and
\cite[Theorem~4.2 and Proposition~4.3]{YuanHeWu}, respectively.
For modules over generalized linear Reedy categories, a lifting of this
kind is proved in \cite[Theorem~9.3]{DiLiLiang} under projectivity
assumptions.
\end{remark}

The following example illustrates Theorem~\ref{thm:directed} with
nonzero multiplication, together with the strictness of the inclusions
\eqref{eq:induced-in-structural}. Here $\operatorname{Add}$ denotes the closure under
arbitrary direct sums and direct summands.

\begin{example}\label{ex:path}
Let $A = k[\varepsilon]/(\varepsilon^2)$ for a field $k$, and let
$Q: 1 \to 2 \to 3$. The forgetful functor
$e: \operatorname{Rep}(Q, A\text{-}\operatorname{Mod}) \to (A\text{-}\operatorname{Mod})^3$ makes
$\operatorname{Rep}(Q, A\text{-}\operatorname{Mod})$ a finite directed cleft extension of
$(A\text{-}\operatorname{Mod})^3$ with
\[
 F(X_1, X_2, X_3) = (0, X_1, X_2 \oplus X_1), \quad
 F^2X = (0, 0, X_1), \quad F^3 = 0,
\]
and with nonzero multiplication $\eta_X: x \mapsto (0, x)$ in the third
coordinate.
For $X = (X_1\xrightarrow aX_2\xrightarrow bX_3)$, the structural complex
has nontrivial parts
\[
 0 \to X_1\xrightarrow aX_2, \qquad
 X_1\xrightarrow{\binom a{-1}}X_2 \oplus X_1
             \xrightarrow{(b, ba)}X_3.
\]
The map $(y, x) \mapsto y+ax$ identifies the quotient of $X_2 \oplus X_1$
by $\operatorname{Im}\binom a{-1}$ with $X_2$. Hence, for
$\mathcal{B} = (A\text{-}\operatorname{Mod})^3$, the class $\Delta(\mathcal{B})$ consists of chains of monomorphisms, while $\nabla(\mathcal{B})$
consists of chains of epimorphisms. Since $A$ is quasi-Frobenius,
Theorem~\ref{thm:directed} recovers the complete hereditary cotorsion triplet
\[
 \bigl(\operatorname{Mon}(Q, A), \operatorname{Rep}(Q, \operatorname{Proj} A), \operatorname{Epi}(Q, A)\bigr)
\]
of \cite[Corollary~6.6]{EPZ}.

Let $S = A/(\varepsilon)$ be the simple $A$-module, let $\sigma: S \to A$ be
given by $\sigma(\bar a) = a\varepsilon$, and let $\pi: A \to S$ be the
quotient map. Then
\[
 M = (S\xrightarrow\sigma A\xrightarrow{\mathrm{id}}A), \qquad
 N = (A\xrightarrow{\mathrm{id}}A\xrightarrow\pi S)
\]
belong to $\Delta(\mathcal{B})$ and $\nabla(\mathcal{B})$, respectively. The exact
sequence $0 \to S\xrightarrow{\sigma}A\xrightarrow\pi S \to 0$ does not split:
$\varepsilon a = 0$ implies $a \in \varepsilon A = \operatorname{Ker}\pi$. On the other hand,
all arrows of induced objects are split monomorphisms, and all arrows of
coinduced objects are split epimorphisms; these properties are preserved by
arbitrary direct sums and direct summands. Thus, for
$\mathcal{B} = (A\text{-}\operatorname{Mod})^3$,
\[
 \operatorname{Add}(l\mathcal{B}) \subsetneq \Delta(\mathcal{B}), \qquad
 \operatorname{Add}(r\mathcal{B}) \subsetneq \nabla(\mathcal{B}),
\]
as witnessed by $M$ and $N$, respectively.
\end{example}

The finite directed condition yields the filtrations used in the proof of
Theorem~\ref{thm:directed}. The following example shows that this lifting can fail under
nilpotence alone, even when $F$ is exact and $F^2 = 0$.

\begin{example}\label{ex:nilpotent-obstruction}
Let $R = k(2\xrightarrow a1)$ and let
$\Lambda = k(2\overset{a, b}{\rightrightarrows}1)$ be the Kronecker algebra.
The projection $\Lambda \to R$ killing $b$ has kernel $J = kb$, with
${}_RJ \cong Re_1$, $J_R \cong e_2R$ and $J \otimes_RJ = 0$. Thus
$\Lambda = R \ltimes J$ and $F = J \otimes_R-$ is exact, with $F^2 = 0$ and
$\eta = 0$.

Consider the complete hereditary cotorsion pair
$(R\text{-}\operatorname{Mod}, \operatorname{Inj} R)$ and the $\Lambda$-module $N$ given by
$N_1 = N_2 = k$ and $a_N = b_N = \mathrm{id}$. Its restriction $eN = (k\xrightarrow{\mathrm{id}}k)$
is injective, and its structure map $F(eN) \to eN$ is monic. Hence
\[
 N \in \Delta(R\text{-}\operatorname{Mod}) \cap e^{-1}(\operatorname{Inj} R).
\]
Nevertheless, $N$ has a nonsplit self-extension: let $E_1 = E_2 = k^2$ and
\[
 a_E = \begin{pmatrix}1 & 0\\0 & 1\end{pmatrix}, \qquad
 b_E = \begin{pmatrix}1 & 1\\0 & 1\end{pmatrix}.
\]
The first coordinates give a submodule $N$ with quotient $N$. A section
with components $(u_1, 1)^t$ and $(u_2, 1)^t$ would satisfy $u_1 = u_2$
from $a_E$ and $u_1 = u_2+1$ from $b_E$. Thus
$\operatorname{Ext}^1_\Lambda(N, N) \neq 0$, so
$(\Delta(R\text{-}\operatorname{Mod}), e^{-1}(\operatorname{Inj} R))$ is not orthogonal and
\[
 {}^\perp e^{-1}(\operatorname{Inj} R) \subsetneq \Delta(R\text{-}\operatorname{Mod}).
\]

Moreover, $N \notin i(\operatorname{Inj} R)$ because $b_N \neq 0$, and its only nonzero
proper submodule restricts to the noninjective simple $R$-module at
vertex $1$. Hence $N$ has no finite filtration with factors in $i(\operatorname{Inj} R)$.
\end{example}

We are now in a position to prove Theorem C.

\begin{theorem}\label{thm:three-levels}
In the setting of Theorem~\ref{thm:directed}, suppose that
\[
 L_1F(\mathcal{C}) = 0, \qquad R^1G(\mathcal{D}) = 0, \qquad
 \mathcal{C}_j \cap \mathcal{Y}_j \subseteq \mathcal{X}_j, \quad \mathcal{X}_j \cap \mathcal{D}_j \subseteq \mathcal{Y}_j
 \quad(1 \leq j \leq n).
\]
Then the following are equivalent:
\begin{enumerate}
 \item $(\mathcal{C}, \mathcal{D})$ is admissible balanced in $\mathcal{B}$;
 \item $(\operatorname{add}(l\mathcal{C}), \operatorname{add}(r\mathcal{D}))$ is admissible balanced in $\mathcal{A}$;
 \item $(\Delta(\mathcal{C}), \nabla(\mathcal{D}))$ is admissible balanced in
 $\mathcal{A}$;
 \item $\mathcal{Y} = \mathcal{X}$;
 \item $\Delta(\mathcal{C})^\perp = {}^\perp\nabla(\mathcal{D})$.
\end{enumerate}
\end{theorem}

\begin{proof}
The equivalence $(1) \Leftrightarrow (2)$ is Theorem~\ref{thm:induced}. Since the product classes satisfy
\eqref{eq:core-compatibility}, the remaining equivalences follow from
Theorems~\ref{thm:directed} and~\ref{thm:structural-criterion}.
\end{proof}

\section{Applications}\label{sec:applications}

Throughout this section, rings are associative with identity, ring homomorphisms preserve identity, and all modules are unital left modules.

\subsection{Cleft extensions of rings}\label{subsec:split}

Let $\alpha:R\to S$ and $\beta:S\to R$ be ring homomorphisms such that $\beta\circ\alpha=\mathrm{id}_R$. We call the triple $(S,\beta,\alpha)$ a \emph{cleft extension} of $R$. Restriction of scalars along $\alpha$ gives a faithful exact functor
\[
e:S\text{-}\operatorname{Mod}\longrightarrow
R\text{-}\operatorname{Mod},
\]
while restriction of scalars along $\beta$ gives an exact functor
$i:R\text{-}\operatorname{Mod}\to S\text{-}\operatorname{Mod}$ with $ei=\mathrm{id}_{R\text{-}\operatorname{Mod}}$.
The induction and coinduction functors
\[
l=S\otimes_R-,
\qquad
r=\operatorname{Hom}_R({}_RS,-)
\]
satisfy $l\dashv e\dashv r$. Thus $(R\text{-}\operatorname{Mod},S\text{-}\operatorname{Mod},i,e,l)$
is a cleft extension of module categories, and the following result follows directly from Theorem~\ref{thm:induced}.
\begin{corollary}\label{cor:split-rings}
	Let $\mathcal{C},\mathcal{D}$ be subcategories of
	$R\text{-}\operatorname{Mod}$.
	Then $(\mathcal{C},\mathcal{D})$ is an admissible balanced pair
	in $R\text{-}\operatorname{Mod}$ if and only if
	\[
	\bigl(
	\operatorname{add}(S\otimes_R\mathcal{C}),
	\operatorname{add}(\operatorname{Hom}_R({}_RS,\mathcal{D}))
	\bigr)
	\]
	is an admissible balanced pair in $S\text{-}\operatorname{Mod}$.
\end{corollary}

Set $J=\operatorname{Ker}\beta$, viewed as an $R$-bimodule via $\alpha$. Multiplication in $S$ induces an associative $R$-bimodule map
$\theta:J\otimes_R J\longrightarrow J$. Let $F=J\otimes_R-$ and $\eta=\theta\otimes_R-$. Then there is an equivalence $S\text{-}\operatorname{Mod}\simeq (R\text{-}\operatorname{Mod})\ltimes_\eta F$, see \cite[Section~3]{BelCleft}.

As an application, consider a Morita context ring
\[
\Lambda=\begin{pmatrix}A&N\\ M&B\end{pmatrix}
\]
whose pairings $N\otimes_B M\to A$ and $M\otimes_A N\to B$ are both zero. The diagonal projection
$\Lambda\to A\times B$ is a ring homomorphism split by the diagonal embedding.

Let $\mathcal{C}_1,\mathcal{D}_1$ be subcategories of $A\text{-}\operatorname{Mod}$
and $\mathcal{C}_2,\mathcal{D}_2$ be subcategories of $B\text{-}\operatorname{Mod}$. Under the canonical equivalence
\[
(A\times B)\text{-}\operatorname{Mod}
\simeq A\text{-}\operatorname{Mod}\times B\text{-}\operatorname{Mod},
\]
approximations and the balance condition are determined componentwise. Thus Corollary~\ref{cor:split-rings} implies that
$(\mathcal{C}_1,\mathcal{D}_1)$ and
$(\mathcal{C}_2,\mathcal{D}_2)$ are admissible balanced pairs in their respective module categories if and only if
\[
\bigl(
\operatorname{add}(l(\mathcal{C}_1\times\mathcal{C}_2)),
\operatorname{add}(r(\mathcal{D}_1\times\mathcal{D}_2))
\bigr)
\]
is an admissible balanced pair in
$\Lambda\text{-}\operatorname{Mod}$. Here $l$ and $r$ are induction and coinduction along the
diagonal embedding $A\times B\to\Lambda$.

When $N=0$, $\Lambda$ is a formal triangular matrix ring.

\subsection{Comma categories and triangular matrix rings}
\label{subsec:triangular}

Let $\mathcal{B}_1$ and $\mathcal{B}_2$ be abelian categories, and let $H:\mathcal{B}_1\to\mathcal{B}_2$ be an additive right exact functor with right adjoint $K$.
The comma category $(H\downarrow\mathcal{B}_2)$ is a cleft extension of $\mathcal{B}_1\times\mathcal{B}_2$, with
\[
e(X,Y,\varphi)=(X,Y),
\qquad
i(X,Y)=(X,Y,0).
\]
The associated endofunctors are
\[
F(X,Y)=(0,HX),
\qquad
G(X,Y)=(KY,0),
\]
so that $F^2=G^2=0$. The induction and coinduction functors are
\begin{equation}\label{eq:comma-adjoints}
	\mathbf p(X,Y)=(X,HX\oplus Y,\kappa),
	\qquad
	\mathbf h(X,Y)=(X\oplus KY,Y,\mathrm{ev}),
\end{equation}
where $\kappa=\binom{\mathrm{id}_{HX}}{0}$ and
$\mathrm{ev}:H(X\oplus KY)\to Y$ is the morphism corresponding under $H\dashv K$ to the projection
$(0,\mathrm{id}_{KY}):X\oplus KY\to KY$.

For a comma object $(X,Y,\varphi)$, write
$\widetilde{\varphi}:X\to KY$ for the adjoint of
$\varphi:HX\to Y$. For subcategories $\mathcal{C}_j,\mathcal{D}_j\subseteq\mathcal{B}_j$,
define
\[
\begin{split}
	\Phi(\mathcal{C}_1,\mathcal{C}_2)
	&=\{(X,Y,\varphi)\in(H\downarrow\mathcal{B}_2)
	\mid \varphi\text{ is a monomorphism},\
	X\in\mathcal{C}_1,\
	\operatorname{Coker}\varphi\in\mathcal{C}_2\},\\
	\Psi(\mathcal{D}_1,\mathcal{D}_2)
	&=\{(X,Y,\varphi)\in(H\downarrow\mathcal{B}_2)
	\mid \widetilde{\varphi}\text{ is an epimorphism},\
	\operatorname{Ker}\widetilde{\varphi}\in\mathcal{D}_1,\
	Y\in\mathcal{D}_2\}.
\end{split}
\]
By \eqref{eq:L1q-model} and \eqref{eq:R1s-model}, we have
\[
\Phi(\mathcal{C}_1,\mathcal{C}_2)
=\Delta(\mathcal{C}_1\times\mathcal{C}_2),
\qquad
\Psi(\mathcal{D}_1,\mathcal{D}_2)
=\nabla(\mathcal{D}_1\times\mathcal{D}_2).
\]

Under our standing convention that subcategories are closed under isomorphisms and direct summands, the classes
$\mathbf p(\mathcal{C}_1,\mathcal{C}_2)$ and
$\mathbf h(\mathcal{D}_1,\mathcal{D}_2)$ are closed under direct summands. Indeed, let $(X',Y',\varphi')$ be a direct summand of
$\mathbf p(X,Y)$, where $X\in\mathcal{C}_1$ and
$Y\in\mathcal{C}_2$. Then $\varphi'$ is a split monomorphism, and hence
\[
(X',Y',\varphi')
\cong\mathbf p(X',\operatorname{Coker}\varphi').
\]
Moreover, $X'$ and $\operatorname{Coker}\varphi'$ are direct summands of $X$ and $\operatorname{Coker}\kappa\cong Y$, respectively. Thus $(X',Y',\varphi')\in\mathbf p(\mathcal{C}_1,\mathcal{C}_2)$. The argument for $\mathbf h$ is dual.

Since the subcategories $\mathcal{C}_j$ and $\mathcal{D}_j$ are closed under finite direct sums, the additivity of $\mathbf p$ and $\mathbf h$ gives
\[
\operatorname{add}\bigl(\mathbf p(\mathcal{C}_1,\mathcal{C}_2)\bigr)
=\mathbf p(\mathcal{C}_1,\mathcal{C}_2),
\qquad
\operatorname{add}\bigl(\mathbf h(\mathcal{D}_1,\mathcal{D}_2)\bigr)
=\mathbf h(\mathcal{D}_1,\mathcal{D}_2).
\]
Since approximations and the balance condition in
$\mathcal{B}_1\times\mathcal{B}_2$ are determined componentwise, Theorem~\ref{thm:induced} yields the following result.
\begin{corollary}\label{cor:comma-induced}
For subcategories $\mathcal{C}_j, \mathcal{D}_j \subseteq \mathcal{B}_j$, $j = 1, 2$, the following are
equivalent:
\begin{enumerate}
 \item $(\mathcal{C}_j, \mathcal{D}_j)$ is an admissible balanced pair in $\mathcal{B}_j$ for
 $j = 1, 2$;
 \item $(\mathbf p(\mathcal{C}_1, \mathcal{C}_2), \mathbf h(\mathcal{D}_1, \mathcal{D}_2))$ is an admissible
 balanced pair in $(H \downarrow \mathcal{B}_2)$.
\end{enumerate}
\end{corollary}

We now apply Theorems~\ref{thm:directed} and~\ref{thm:three-levels} to the comma category.
\begin{corollary}\label{cor:comma}
Suppose that $\mathcal{B}_1$ and $\mathcal{B}_2$ have enough projective and injective objects. For $j = 1, 2$, let $(\mathcal{C}_j, \mathcal{Y}_j)$ and $(\mathcal{X}_j, \mathcal{D}_j)$ be complete
hereditary cotorsion pairs in $\mathcal{B}_j$. If
\begin{equation}\label{eq:comma-relative}
 L_1H(\mathcal{C}_1) = 0, \qquad R^1K(\mathcal{D}_2) = 0,
\end{equation}
then
\[
 \bigl(\Phi(\mathcal{C}_1, \mathcal{C}_2), e^{-1}(\mathcal{Y}_1 \times \mathcal{Y}_2)\bigr), \qquad
 \bigl(e^{-1}(\mathcal{X}_1 \times \mathcal{X}_2), \Psi(\mathcal{D}_1, \mathcal{D}_2)\bigr)
\]
are complete hereditary cotorsion pairs in $(H\downarrow\mathcal{B}_2)$. If, in addition,
$\mathcal{C}_j \cap \mathcal{Y}_j \subseteq \mathcal{X}_j$ and $\mathcal{X}_j \cap \mathcal{D}_j \subseteq \mathcal{Y}_j$ for $j = 1, 2$, then the following are equivalent:
\begin{enumerate}
 \item $(\mathcal{C}_j, \mathcal{D}_j)$ is admissible balanced in $\mathcal{B}_j$ for $j = 1, 2$;
 \item $(\mathbf p(\mathcal{C}_1, \mathcal{C}_2), \mathbf h(\mathcal{D}_1, \mathcal{D}_2))$ is
 admissible balanced;
 \item $(\Phi(\mathcal{C}_1, \mathcal{C}_2), \Psi(\mathcal{D}_1, \mathcal{D}_2))$ is admissible balanced;
 \item $\mathcal{Y}_j = \mathcal{X}_j$ for $j = 1, 2$;
 \item $\Phi(\mathcal{C}_1, \mathcal{C}_2)^\perp = {}^\perp\Psi(\mathcal{D}_1, \mathcal{D}_2)$.
\end{enumerate}
\end{corollary}

\begin{proof}
	The comma category is a finite directed extension of $\mathcal{B}_1\times\mathcal{B}_2$, with $F_{21}=H$ and $G_{12}=K$. Since
	\[
	L_1F(X,Y)\cong(0,L_1H(X)),\qquad
	R^1G(X,Y)\cong(R^1K(Y),0),
	\]
	condition~\eqref{eq:comma-relative} gives
	$L_1F(\mathcal{C}_1\times\mathcal{C}_2)=0$ and
	$R^1G(\mathcal{D}_1\times\mathcal{D}_2)=0$.
	Thus Theorem~\ref{thm:directed}, together with the identities
	\[
	\Delta(\mathcal{C}_1\times\mathcal{C}_2)
	=\Phi(\mathcal{C}_1,\mathcal{C}_2),\qquad
	\nabla(\mathcal{D}_1\times\mathcal{D}_2)
	=\Psi(\mathcal{D}_1,\mathcal{D}_2),
	\]
	yields the asserted complete hereditary cotorsion pairs.
	
	Under the additional inclusions, 
	Theorem~\ref{thm:three-levels} gives the equivalence of \textup{(1)--(5)}.
\end{proof}

The cotorsion-pair assertions of Corollary~\ref{cor:comma} are those of
\cite{HuZhuComma} and \cite{YuanHeWu} recalled in
Remark~\ref{rem:hu-comparison}.

In particular, taking $N = 0$ and writing $U = M$ in the Morita context ring
of Subsection~\ref{subsec:split} gives the triangular matrix ring
$T = \left(\begin{smallmatrix}A & 0\\ U & B\end{smallmatrix}\right)$. Its
module category $T\text{-}\operatorname{Mod}$ is the comma category for $H = U \otimes_A-$, with
right adjoint $K = \operatorname{Hom}_B(U, -)$, and \eqref{eq:comma-relative} becomes
\begin{equation}\label{eq:triangular-relative}
 \operatorname{Tor}^A_1(U, \mathcal{C}_1) = 0, \qquad \operatorname{Ext}^1_B(U, \mathcal{D}_2) = 0.
\end{equation}
Since $\mathcal{C}_1$ is resolving and $\mathcal{D}_2$ is coresolving, dimension shifting
gives $\operatorname{Tor}_m^A(U, \mathcal{C}_1) = 0$ and $\operatorname{Ext}_B^m(U, \mathcal{D}_2) = 0$ for all $m \geq 1$, which
are the vanishing conditions used in \cite{Mao}. In this case
Corollary~\ref{cor:comma-induced} recovers \cite[Theorem~2.2]{Mao} for the
additive, summand-closed classes considered here, and
Corollary~\ref{cor:comma} recovers \cite[Theorem~2.5]{Mao}.
We next consider the arrow category
$\operatorname{Arr}(\mathcal{B})$, whose objects are morphisms $X\to Y$ in $\mathcal{B}$. This is the comma category associated with
$\mathrm{id}_{\mathcal{B}}$.

\begin{corollary}\label{cor:arrows}
Let $\mathcal{B}$ have enough projective and injective objects, and suppose that
$(\mathcal{C}, \mathcal{C}^\perp)$ and $({}^\perp\mathcal{D}, \mathcal{D})$ are complete hereditary cotorsion
pairs. The following are equivalent:
\begin{enumerate}
 \item $(\Phi(\mathcal{C}, \mathcal{C}), \Psi(\mathcal{D}, \mathcal{D}))$ is admissible balanced in
 $\operatorname{Arr}(\mathcal{B})$;
 \item $\mathcal{C}^\perp = {}^\perp\mathcal{D}$.
\end{enumerate}
In this case $(\mathcal{C}, \mathcal{C}^\perp, \mathcal{D})$ is a complete hereditary cotorsion
triplet and $(\mathcal{C}, \mathcal{D})$ is admissible balanced.
\end{corollary}

\begin{proof}
	Suppose that \textup{(2)} holds. Since the identity functor is
	exact, the cotorsion-pair assertions of
	Corollary~\ref{cor:comma} give a complete hereditary cotorsion triplet
	\[
	\bigl(
	\Phi(\mathcal{C},\mathcal{C}),
	e^{-1}(\mathcal{C}^{\perp}\times\mathcal{C}^{\perp}),
	\Psi(\mathcal{D},\mathcal{D})
	\bigr)
	\]
	in $\operatorname{Arr}(\mathcal{B})$.
	Thus \textup{(1)} follows from
	Lemma~\ref{lem:triplet-balanced}.
	
	Conversely, suppose that \textup{(1)} holds.
	Let  $M\in\mathcal{B}$. Consider the short exact sequence
	\begin{equation}\label{eq:arrow-test}
		0\longrightarrow (0\to M)
		\longrightarrow (M\xrightarrow{\mathrm{id}_M}M)
		\longrightarrow (M\to0)
		\longrightarrow0.
	\end{equation}
	For a monomorphism $u:U\to V$, applying
	$\operatorname{Hom}_{\operatorname{Arr}(\mathcal{B})}
	((U\xrightarrow{u}V),-)$ to \eqref{eq:arrow-test} gives the map
	\[
	\operatorname{Hom}_{\mathcal{B}}(V,M)
	\xrightarrow{u^*}
	\operatorname{Hom}_{\mathcal{B}}(U,M).
	\]
	Dually, for an epimorphism $v:V\to W$, applying
	$\operatorname{Hom}_{\operatorname{Arr}(\mathcal{B})}
	(-,(V\xrightarrow{v}W))$ gives the map
	\[
	\operatorname{Hom}_{\mathcal{B}}(M,V)
	\xrightarrow{v_*}
	\operatorname{Hom}_{\mathcal{B}}(M,W).
	\]
	The sequence \eqref{eq:arrow-test} is
	$\operatorname{Hom}((U\xrightarrow{u}V),-)$-exact if and only if $u^*$ is surjective, and it is
	$\operatorname{Hom}(-,(V\xrightarrow{v}W))$-exact if and only if $v_*$ is surjective.
	
	Let $M\in{}^\perp\mathcal{D}$. For every $(V\xrightarrow{v}W)\in \Psi(\mathcal{D},\mathcal{D})$, we have
	$\operatorname{Ker}v\in\mathcal{D}$, so $v_*$ is surjective. Hence \eqref{eq:arrow-test} is
	$\operatorname{Hom}(-,\Psi(\mathcal{D},\mathcal{D}))$-exact. By balance, it is also
	$\operatorname{Hom}(\Phi(\mathcal{C},\mathcal{C}),-)$-exact.
	
	For $C\in\mathcal{C}$, choose a short exact sequence
	\[
	0\longrightarrow L\xrightarrow{u}P
	\longrightarrow C\longrightarrow0
	\]
	with $P$ projective.
	Since $\mathcal{C}$ is resolving, $L\in\mathcal{C}$, and
	therefore $(L\xrightarrow{u}P)\in
	\Phi(\mathcal{C},\mathcal{C})$.
	Consequently, $u^*$ is surjective, and
	\[
	\operatorname{Ext}^1_{\mathcal{B}}(C,M)
	\cong\operatorname{Coker}u^*=0.
	\]
	Thus ${}^\perp\mathcal{D}\subseteq\mathcal{C}^\perp$.
	
	Now let $M\in\mathcal{C}^\perp$. For every $(U\xrightarrow{u}V)\in \Phi(\mathcal{C},\mathcal{C})$, we have
	$\operatorname{Coker}u\in\mathcal{C}$, so $u^*$ is surjective. Thus \eqref{eq:arrow-test} is
	$\operatorname{Hom}(\Phi(\mathcal{C},\mathcal{C}),-)$-exact and, by balance,
	$\operatorname{Hom}(-,\Psi(\mathcal{D},\mathcal{D}))$-exact.
	
	For $D\in\mathcal{D}$, choose a short exact sequence
	\[
	0\longrightarrow D\longrightarrow I
	\xrightarrow{v}L\longrightarrow0
	\]
	with $I$ injective. Since $\mathcal{D}$ is coresolving, $L\in\mathcal{D}$, and
	therefore $(I\xrightarrow{v}L)\in
	\Psi(\mathcal{D},\mathcal{D})$. Consequently, $v_*$ is surjective, and
	\[
	\operatorname{Ext}^1_{\mathcal{B}}(M,D)
	\cong\operatorname{Coker}v_*=0.
	\]
	Hence $\mathcal{C}^\perp\subseteq{}^\perp\mathcal{D}$, proving \textup{(2)}.
	
	Finally, the equality
	$\mathcal{C}^\perp={}^\perp\mathcal{D}$
	and the hypotheses give the complete hereditary cotorsion triplet $(\mathcal{C},\mathcal{C}^\perp,\mathcal{D})$.
	Thus $(\mathcal{C},\mathcal{D})$ is an admissible balanced pair by Lemma~\ref{lem:triplet-balanced}.
\end{proof}

For $\mathcal{B}=R\text{-}\operatorname{Mod}$,
Corollary~\ref{cor:arrows} recovers the equivalence of \textup{(2)} and \textup{(3)} in
\cite[Corollary~2.6]{Mao}. It extends this equivalence to abelian categories with enough projective and injective objects.

\subsection{Gorenstein specializations}\label{subsec:gorenstein}

Write $\operatorname{GProj} R$ and $\operatorname{GInj} R$ for the classes of Gorenstein projective
and Gorenstein injective left $R$-modules. Keep the triangular matrix
ring $T$ of Subsection~\ref{subsec:triangular}, and suppose that $U_A$
has finite flat dimension and ${}_BU$ has finite projective dimension.
Then
\[
 \operatorname{GProj} T = \Phi(\operatorname{GProj} A, \operatorname{GProj} B), \qquad
 \operatorname{GInj} T = \Psi(\operatorname{GInj} A, \operatorname{GInj} B).
\]
The first equality is \cite[Theorem~1.1]{LiZhengHuZhu}, stated there for
the upper triangular matrix ring
$\left(\begin{smallmatrix}B & U\\ 0 & A\end{smallmatrix}\right) \cong T$, and the
second is \cite[Propositions~3.3, 5.2 and~5.5]{YuanHeWu}.

Recall that a ring $R$ is left virtually Gorenstein if
\[
 (\operatorname{GProj} R)^\perp = {}^\perp(\operatorname{GInj} R)
\]
\cite[Section~4]{EPZ}. By
\cite[Corollary~3.11]{EHT}, applied to the opposite ring, this is equivalent
to balance of $(\operatorname{GProj} R, \operatorname{GInj} R)$ and yields a complete hereditary
cotorsion triplet. We obtain the following characterization for $T$.

\begin{corollary}\label{cor:virtually-gorenstein}
Let $A$ and $B$ be rings, and suppose that $U_A$ has finite flat
dimension and ${}_BU$ has finite projective dimension. The following are
equivalent:
\begin{enumerate}
 \item $(\operatorname{GProj} A, \operatorname{GInj} A)$ and $(\operatorname{GProj} B, \operatorname{GInj} B)$ are admissible
 balanced;
 \item $(\operatorname{GProj} T, \operatorname{GInj} T)$ is admissible balanced;
 \item $(\mathbf p(\operatorname{GProj} A, \operatorname{GProj} B), \mathbf h(\operatorname{GInj} A, \operatorname{GInj} B))$ is
 admissible balanced;
 \item $A$ and $B$ are left virtually Gorenstein;
 \item $T$ is left virtually Gorenstein.
\end{enumerate}
In particular, Theorem D holds.
\end{corollary}

\begin{proof}
Since $\operatorname{GProj} R\supseteq\operatorname{Proj} R$ and $\operatorname{GInj} R\supseteq\operatorname{Inj} R$, every
$\operatorname{GProj} R$-precover is epic and every $\operatorname{GInj} R$-preenvelope is monic.
Thus balance in the sense of \cite[Definition~3.1]{EHT} is admissible
balance in the sense of Lemma~\ref{lem:balance-criterion}.
Corollary~\ref{cor:comma-induced} gives
$(1) \Leftrightarrow (3)$.  Applying \cite[Corollary~3.11]{EHT} to the opposite rings $A^{\mathrm{op}}$, $B^{\mathrm{op}}$ and $T^{\mathrm{op}}$ gives $\textup{(1)}\Leftrightarrow\textup{(4)}$ and
$\textup{(2)}\Leftrightarrow\textup{(5)}$. It remains to compare \textup{(1)} and \textup{(2)}.

Choose nonnegative integers $d, p$ with $\operatorname{fd}(U_A) \leq d$ and
$\operatorname{pd}({}_BU) \leq p$. For $C \in \operatorname{GProj} A$ and $D \in \operatorname{GInj} B$,
complete projective resolutions and complete injective coresolutions
give exact sequences $0 \to C \to P \to C^1 \to 0$ and
$0 \to D^{-1} \to I \to D \to 0$ in which all successive cycles are again
Gorenstein projective and Gorenstein injective, respectively. Dimension
shifting gives
\[
 \operatorname{Tor}^A_1(U, C) \cong \operatorname{Tor}^A_{d+1}(U, C^d) = 0, \qquad
 \operatorname{Ext}^1_B(U, D) \cong \operatorname{Ext}^{p+1}_B(U, D^{-p}) = 0,
\]
so that \eqref{eq:triangular-relative} holds. If (1) holds,
\cite[Corollary~3.11]{EHT} supplies complete hereditary cotorsion
triplets
\[
 (\operatorname{GProj} R, \mathcal{W}_R, \operatorname{GInj} R), \qquad
 \mathcal{W}_R = (\operatorname{GProj} R)^\perp = {}^\perp(\operatorname{GInj} R)
 \quad(R = A, B).
\]
By Corollary~\ref{cor:comma} and the descriptions of $\operatorname{GProj} T$ and $\operatorname{GInj} T$ above,
\[
 \bigl(\operatorname{GProj} T, e^{-1}(\mathcal{W}_A \times \mathcal{W}_B), \operatorname{GInj} T\bigr)
\]
is a complete hereditary cotorsion triplet, and
Lemma~\ref{lem:triplet-balanced} proves (2). Conversely,
$(2) \Rightarrow (1)$ is Lemma~\ref{lem:reflection}.
\end{proof}

For Artin algebras, the equivalence of \textup{(4)} and \textup{(5)} is
\cite[Proposition~3.3]{EHSL}; the four conditions imposed there on $U$
follow from our assumptions by dimension shifting. The implication
$\textup{(2)}\Rightarrow\textup{(1)}$ is \cite[Corollary~2.7]{Mao}, and
Corollary~\ref{cor:virtually-gorenstein} recovers \cite[Theorem~2.8]{Mao}
when $U_A$ is flat, ${}_BU$ is projective, and $A$ and $B$ are left
noetherian rings of finite left self-injective dimension; none of these
assumptions is needed here.

\begin{corollary}\label{cor:lower-triangular}
Let $R$ be a ring and $n \geq 1$, and let $T_n(R)$ be the ring of lower
triangular $n \times n$ matrices over $R$. Then $T_n(R)$ is left virtually
Gorenstein if and only if $R$ is left virtually Gorenstein.
\end{corollary}

\begin{proof}
The case $n = 1$ is trivial. For $n \geq 2$,
\[
 T_n(R) = \begin{pmatrix}T_{n-1}(R) & 0\\ U & R\end{pmatrix},
\]
where $U = R^{n-1}$ is the last row of $T_n(R)$ with its diagonal entry
removed. As a left $R$-module $U$ is free, and as a right
$T_{n-1}(R)$-module it is isomorphic to the last row of $T_{n-1}(R)$,
which is a direct summand of $T_{n-1}(R)$. Thus ${}_RU$ and
$U_{T_{n-1}(R)}$ are projective, and
Corollary~\ref{cor:virtually-gorenstein} shows that $T_n(R)$ is left
virtually Gorenstein if and only if $T_{n-1}(R)$ and $R$ are left
virtually Gorenstein. The assertion follows by induction on $n$.
\end{proof}

For $n=2$, this recovers the equivalence of
\textup{(4)} and \textup{(5)} in
\cite[Corollary~2.9]{Mao}, without assuming that $R$ is left noetherian or has finite left self-injective dimension.
For $n = 2$ and $R$ an Artin algebra, it is \cite[Example~3.4]{EHSL}.

\section*{Declaration of AI Use}

Claude Fable 5.2 (Anthropic) was used in the final stage of manuscript preparation to check notation, suggest wording improvements, and help identify possible gaps in proofs and missing references. All results and proofs are due to the authors, who independently checked every argument and take full responsibility for the manuscript.

\section*{Acknowledgements}

The first and third authors were supported by the China Scholarship
Council (Grant Nos.\ 202606620154 and 202606620170, respectively). The
first author is grateful to the University of Murcia for its hospitality
during the preparation of this article.

\bigskip
{\small\raggedright
\noindent\textbf{Junpeng Ren}\\
School of Mathematics and Statistics, Northeast Normal University\\
Changchun 130024, China\\
E-mail: renjp@nenu.edu.cn

\medskip
\noindent\textbf{Zheng-Gang Wu}\\
School of Mathematical Sciences, Beijing Normal University\\
Beijing 100875, China\\
E-mail: 202631130011@mail.bnu.edu.cn

\medskip
\noindent\textbf{Mengzhen Zhu}\\
School of Mathematics and Statistics, Northeast Normal University\\
Changchun 130024, China\\
E-mail: mzzhu@nenu.edu.cn
\par}


\begin{thebibliography}{99}

\bibitem{BelCleft}
A.~Beligiannis,
\emph{Cleft extensions of abelian categories and applications to ring theory},
Comm. Algebra \textbf{28} (2000), 4503--4546.

\bibitem{BelRelative}
A.~Beligiannis,
\emph{On the relative homology of cleft extensions of rings and abelian categories},
J. Pure Appl. Algebra \textbf{150} (2000), 237--299.

\bibitem{BeligiannisReiten}
A.~Beligiannis and I.~Reiten,
\emph{Homological and homotopical aspects of torsion theories},
Mem. Amer. Math. Soc. \textbf{188} (2007), no.~883.

\bibitem{Chen}
X.-W.~Chen,
\emph{Homotopy equivalences induced by balanced pairs},
J. Algebra \textbf{324} (2010), 2718--2731.



\bibitem{ChenQinYang}
P.~Chen, Y.~Qin and J.~Yang,
\emph{Cotorsion pairs in $\eta$-extensions with applications to $\theta$-extension rings and tensor rings},
\href{https://arxiv.org/abs/2611.19449v1}{arXiv:2611.19449v1} (2026).

\bibitem{DaleziosStovicek}
G.~Dalezios and J.~\v{S}\v{t}ov\'i\v{c}ek,
\emph{Linear Reedy categories, quasi-hereditary algebras and model structures},
Adv. Math. \textbf{481} (2025), Paper No.~110550.

\bibitem{DiLiLiang}
Z.~Di, L.~Li and L.~Liang,
\emph{Representations of generalized linear Reedy categories and abelian model structures},
\href{https://arxiv.org/abs/2601.01187v2}{arXiv:2601.01187v2} (2026).

\bibitem{EJTX}
E.~E.~Enochs, O.~M.~G.~Jenda, B.~Torrecillas and J.~Z.~Xu,
\emph{Torsion theory relative to Ext},
Technical Report 98-11, Department of Mathematics, University of Kentucky, 1998.

\bibitem{EnochsJenda}
E.~E.~Enochs and O.~M.~G.~Jenda,
\emph{Relative Homological Algebra},
de Gruyter Expositions in Mathematics, vol.~30,
Walter de Gruyter, Berlin, 2000.

\bibitem{EHSL}
H.~Eshraghi, R.~Hafezi, Sh.~Salarian and Z.~W.~Li,
\emph{Gorenstein projective modules over triangular matrix rings},
Algebra Colloq. \textbf{23} (2016), 97--104.

\bibitem{EHT}
S.~Estrada, J.~Hu and J.~Trlifaj,
\emph{Balanced pairs, virtually Gorenstein rings, and cotorsion torsion triples},
\href{https://arxiv.org/abs/2602.18818v2}{arXiv:2602.18818v2} (2026).

\bibitem{EPZ}
S.~Estrada, M.~A.~P\'erez and H.~Zhu,
\emph{Balanced pairs, cotorsion triplets and quiver representations},
Proc. Edinb. Math. Soc. (2) \textbf{63} (2020), 67--90.

\bibitem{FGR}
R.~M.~Fossum, P.~A.~Griffith and I.~Reiten,
\emph{Trivial Extensions of Abelian Categories},
Lecture Notes in Mathematics, vol.~456, Springer-Verlag, Berlin, 1975.

\bibitem{FuHuZhangZhu}
X.~Fu, J.~Hu, D.~Zhang and H.~Zhu,
\emph{Balanced pairs on triangulated categories},
Algebra Colloq. \textbf{30} (2023), 385--394.

\bibitem{GreenPsaroudakisSolberg}
E.~L.~Green, C.~Psaroudakis and {\O}.~Solberg,
\emph{Reduction techniques for the finitistic dimension},
Trans. Amer. Math. Soc. \textbf{374} (2021), 6839--6879.

\bibitem{HeZhou}
J.~He and P.~Zhou,
\emph{Balanced pairs and recollements in extriangulated categories with negative first extensions},
J. Algebra Appl. \textbf{24} (2025), Paper No.~2550110.

\bibitem{HolmJorgensen}
H.~Holm and P.~J{\o}rgensen,
\emph{Cotorsion pairs in categories of quiver representations},
Kyoto J. Math. \textbf{59} (2019), 575--606.

\bibitem{Hu}
D.~Hu,
\emph{Cotorsion pairs in extensions of abelian categories},
Bull. Malays. Math. Sci. Soc. \textbf{49} (2026), Paper No.~58.

\bibitem{HuZhuComma}
J.~Hu and H.~Zhu,
\emph{Special precovering classes in comma categories},
Sci. China Math. \textbf{65} (2022), 933--950.

\bibitem{Kostas}
P.~Kostas,
\emph{Cleft extensions of rings and singularity categories},
J. Algebra \textbf{685} (2026), 160--224.

\bibitem{LiZhengHuZhu}
H.~Li, Y.~Zheng, J.~Hu and H.~Zhu,
\emph{Gorenstein projective modules and recollements over triangular matrix rings},
Comm. Algebra \textbf{48} (2020), 4932--4947.

\bibitem{LiangMaYang}
L.~Liang, Y.~Ma and G.~Yang,
\emph{Homological invariant properties under cleft extensions},
\href{https://arxiv.org/abs/2506.02691v2}{arXiv:2506.02691v2} (2025).

\bibitem{LiuGao}
D.~Liu and H.~Gao,
\emph{Higher cotorsion pairs on image-nilpotent $\eta$-extensions of abelian categories},
\href{https://arxiv.org/abs/2609.39003v1}{arXiv:\allowbreak2609.39003v1} (2026).

\bibitem{MaZhengLiu}
Y.~Ma, J.~Zheng and Y.-Z.~Liu,
\emph{Three homological invariants under cleft extensions},
J. Algebra \textbf{694} (2026), 287--323.

\bibitem{Mao}
L.~Mao,
\emph{Admissible balanced pairs over formal triangular matrix rings},
Bull. Korean Math. Soc. \textbf{58} (2021), 1387--1400.

\bibitem{Marmaridis}
N.~Marmaridis,
\emph{On extensions of abelian categories with applications to ring theory},
J. Algebra \textbf{156} (1993), 50--64.

\bibitem{Qin}
Y.~Qin,
\emph{Gorenstein projective objects over cleft extensions},
\href{https://arxiv.org/abs/2507.09109v1}{arXiv:2507.09109v1} (2025).

\bibitem{WangLiHuang}
J.~Wang, H.~Li and Z.~Huang,
\emph{Applications of exact structures in abelian categories},
Publ. Math. Debrecen \textbf{88} (2016), 269--286.

\bibitem{XuFu}
R.~Xu and X.~Fu,
\emph{Ideal balanced pairs in extriangulated categories},
Sci. China Math. \textbf{68} (2025), 271--284.

\bibitem{YuanHeWu}
Y.~Yuan, J.~He and D.~Wu,
\emph{Cotorsion pairs in comma categories},
Czechoslovak Math. J. \textbf{74(149)} (2024), 715--734.

\bibitem{ZhaoSunSilting}
G.-Q.~Zhao and J.-X.~Sun,
\emph{Silting and tilting objects in cleft extensions of abelian categories},
\href{https://arxiv.org/abs/2602.07368v1}{arXiv:2602.07368v1} (2026).

\end{thebibliography}
\end{document}